\documentclass[twoside]{article}

\usepackage[accepted]{aistats2020}
\usepackage[round]{natbib}
\usepackage{amsmath, amsthm, amssymb, amsfonts} 
\usepackage{mathtools}
\usepackage{MnSymbol}
\usepackage{bm}
\usepackage[USenglish]{babel}

\usepackage[ruled, lined]{algorithm2e}
\usepackage{tabularx}
\usepackage{placeins}
\usepackage{caption}
\usepackage{nicefrac}

\newcommand{\commentout}[1]{}

\newcommand{\diag}{{\text{diag}}}

\newcommand{\PP}{\mathbb{P}}
\newcommand{\EE}{{\mathbb{E}}}

\newcommand{\eps}{\varepsilon}

\newcommand{\RR}{\mathbb{R}}

\newcommand{\eqd}{\overset{\mathrm{d}}{=}}

\newcommand{\given}{\, |\, }

\newtheorem{thm}{Theorem}[section]

\DeclareMathOperator{\ex}{\textbf{E}}

\newcommand\mydots{\hbox to 1em{.\hss.\hss.}}

\newcommand{\loss}{\ensuremath{\text{\sc loss}}}

\newcommand{\cB}{\ensuremath{\mathcal B}}

\newcommand{\cF}{\ensuremath{\mathcal F}}

\newcommand{\cN}{\ensuremath{\mathcal N}}

\newcommand{\cQ}{\ensuremath{\mathcal Q}}

\newcommand{\cX}{\ensuremath{\mathcal X}}
\newcommand{\cY}{\ensuremath{\mathcal Y}}
\newcommand{\cZ}{\ensuremath{\mathcal Z}}

\newcommand{\Exp}{\ensuremath{\text{\rm E}}}

\newcommand{\reals}{{\mathbb R}}
\newcommand{\naturals}{{\mathbb N}}

\newcommand{\dhel}{{d}}
\newcommand{\ind}[1]{\mathop{{\bf 1}_{\{#1\}}}}
\newcommand{\Pt}{P} 
\newcommand{\E}{\operatorname{\mathbf{E}}}

\renewcommand\Exp{\ensuremath{\mathbf E}}

\newcommand{\rv}[1]{\underline{#1}}

\newcommand{\fstar}{f^*}
\newcommand{\fopt}{\fstar}

\renewcommand{\loss}{\ensuremath{\ell}}

\newcommand{\xslosslong}[1]{\loss_{#1} - \loss_{\fopt}}
\newcommand{\xsloss}[1]{L_{#1}}
\newcommand{\xslossat}[2]{L_{#1}(#2)}
\newcommand{\xsrisk}[1]{\E [ \xsloss{#1} ]} 
\newcommand{\xsrisklong}[1]{\E [ \xslosslong{#1} ]} 
\newcommand{\xslosstheta}{\loss_\theta - \loss_{\theta^*}}

\newcommand{\dol}{\ensuremath{\Pi}} 
\newcommand{\dolest}{\dol_|} 
\newcommand{\Prior}{\ensuremath{\Pi_0}} 
\newcommand{\KL}{\text{\sc KL}}
\newcommand{\smtuple}{(\Pt,\loss,\model)}

\newcommand{\model}{\mathcal{F}}

\newcommand{\rsc}{\mathrm{IC}} 

\DeclareMathOperator*{\pipes}{\|}
\DeclareMathOperator*{\opmax}{\vee}

\DeclareMathOperator*{\argmin}{arg\,min}

\newtheorem{theorem}{Theorem}

\newtheorem{lemma}{Lemma}
\newtheorem{proposition}{Proposition}
\newtheorem{definition}{Definition}

\newtheorem{example}{Example}

\makeatletter
\DeclareRobustCommand{\qed}{%
	\ifmmode 
	\else \leavevmode\unskip\penalty9999 \hbox{}\nobreak\hfill
	\fi
	\quad\hbox{\qedsymbol}}

\makeatother

\DeclareRobustCommand{\VANDER}[3]{#2}

\usepackage{todonotes}
\usepackage[draft,multiuser]{fixme}
\fxusetheme{colorsig}
\FXRegisterAuthor{peter}{Peter}{\color{blue}PG} 
\FXRegisterAuthor{niche}{Niche}{\color{red}NM}

\begin{document}

%

%

\twocolumn[

\aistatstitle{Safe-Bayesian Generalized Linear Regression}

\aistatsauthor{ Rianne de Heide \And Alisa Kirichenko \And  Nishant A.\ Mehta \And Peter D.\ Gr{\"u}nwald }
\aistatsaddress{ Leiden University \& CWI \And  University of Oxford \And University of Victoria \And CWI \& Leiden University } 
]

\begin{abstract}
  We study generalized Bayesian inference under misspecification,
  i.e.\ when the model is `wrong but useful'. Generalized Bayes equips
  the likelihood with a learning rate $\eta$. We show that for
  generalized linear models (GLMs), $\eta$-generalized Bayes
  concentrates around the best approximation of the truth within the
  model for specific $\eta \neq 1$, even under severely misspecified
  noise, as long as the tails of the true distribution are exponential. We
  derive MCMC samplers for generalized Bayesian lasso and
  logistic regression and give examples of both
  simulated and real-world data in which generalized Bayes
  substantially outperforms standard Bayes.
\end{abstract}
	
\section{INTRODUCTION}
Over the last ten years it has become abundantly clear that Bayesian
inference can behave quite badly under misspecification, i.e., if the
model $\cF$ under consideration is `wrong but useful'
\citep{GrunwaldL07,ErvenGR07,muller2013risk,syring2017calibrating,yao2018using,holmes2017assigning,grunwald2017inconsistency}. For
example, \citet{GrunwaldL07} exhibit a simple nonparametric
classification setting in which, even though the prior puts positive
mass on the unique distribution in $\cF$ that is closest in KL
divergence to the data generating distribution $P$, the posterior
never concentrates around this distribution.
\citet{grunwald2017inconsistency} give a simple misspecified setting
in which standard Bayesian ridge regression, model selection
and model averaging severely overfit small-sample data. 

\citet{grunwald2017inconsistency}
also propose a
remedy for this problem: equip the likelihood with an exponent or {\em
  learning rate\/} $\eta$ (see~\eqref{eq:bayesgenpost} below). 
Such a
{\em generalized Bayesian\/} (also known as {\em fractional\/} or {\em
  tempered\/} Bayesian) approach was considered earlier by
e.g.~\cite{barron1991minimum,WalkerH02,zhang2006information}. In
practice, $\eta$ will usually (but not always --- see
Section~\ref{sec:lassoandhs}
below) be chosen smaller than one, making the prior have a stronger
regularizing influence. \citet{grunwald2017inconsistency} show that for
Bayesian ridge regression and model selection/averaging, this results
in excellent performance, being competitive with standard Bayes if the
model is correct and very significantly outperforming standard Bayes
if it is not. 
Extending Zhang's (\citeyear{zhang2006epsilon,zhang2006information}) earlier work, \citet{grunwald2016fast} (GM from now on) show that, under what was earlier called the $\bar\eta$-{\em central condition\/} (Definition~\ref{def:central} below), generalized Bayes with a specific finite learning rate $\bar\eta$ (usually $\neq 1$) will indeed concentrate in the neighborhood of the `best' $f \in \cF$ with high probability. Here, the `best' $f$  means the one closest in KL divergence to $P$.

Yet, three
important parts of the story are missing in this existing 
work: (1) Can Gr\"unwald-Van Ommen-type examples, showing failure of
standard Bayes $(\eta = 1)$ and empirical success of generalized Bayes
with the right $\eta$, be given more generally, for different priors
$\pi$ (say of lasso-type ($\pi(f) \propto \exp(- \lambda \| f\|_1)$)
rather than ridge-type $\pi(f) \propto \exp(-\lambda \| f \|_2^2)$),
and for different models, say for {\em generalized\/} linear models
(GLMs)? (2) Can we find examples of generalized Bayes outperforming
standard Bayes with real-world data rather than with toy problems such
as those considered by Gr\"unwald and Van Ommen? (3) Does the central
condition --- which allows for good theoretical behavior of
generalized Bayes --- hold for GLMs, under reasonable further conditions?

We answer all three questions in the affirmative: in
Section~\ref{sec:bad} below, we give (a) a toy example on which the
Bayesian lasso and the Horseshoe estimator fail; later in the paper,
in Section~\ref{sec:examples} we also (b) give a toy example on which
standard Bayes logistic regression fails, and (c) two real-world data
sets on which Bayesian lasso and Horseshoe regression fail; in all
cases, (d) generalized Bayes with the right $\eta$ shows much better
performance. In Section~\ref{sec:genglm}, we show (e) that for GLMs,
even if the noise is severely misspecified, as long as the
distribution of the predictor variable $Y$ has exponentially small
tails (which is automatically the case in classification, where the
domain of $Y$ is finite), the central condition holds for some
$\eta > 0$. In combination with (e), GM's existing theoretical results
suggest that generalized Bayes with this $\eta$ should lead to good
results --- this is corroborated by our experimental results in
Section~\ref{sec:examples}.  These findings are not obvious: one might
for example think that the sparsity-inducing prior used by Bayesian
lasso regression circumvents the need for the additional
regularization induced by taking an $\eta < 1$, especially since in the original
setting of Gr\"unwald and Van Ommen, the standard Bayesian lasso
$(\eta =1)$ succeeds. Yet, Example~\ref{ex:basic} below shows that
under a modification of their example, it fails after all. In order to
demonstrate the failure of standard Bayes and the success of
generalized Bayes, we devise (in Section~\ref{sec:sampling}) MCMC
algorithms (f) for generalized Bayes posterior sampling for Bayesian
lasso and logistic regression. (a)-(f) are all novel contributions.

In Section~\ref{sec:setting} we first  define our setting more
precisely. Section~\ref{sec:bad}) gives a first example of bad
standard-Bayesian behavior and Section~\ref{sec:consistency}) recalls a
theorem from GM indicating that under the $\bar\eta$-central
condition, generalized Bayes for $\eta < \bar\eta$ should perform
well. We present our new theoretical results in Section~\ref{sec:genglm}.  We next
(Section~\ref{sec:sampling}), present our algorithms for generalized
Bayesian posterior sampling, and we continue
(Section~\ref{sec:examples}) to empirically demonstrate how generalized
Bayes outperforms standard Bayes under misspecification. All proofs are in Appendix~\ref{app:proofs}.
\section{THE SETTING}\label{sec:setting}
A \emph{learning problem} can
be characterized by a tuple $(P, \loss, \cF)$, where $\cF$ is a set of
predictors, also referred to as a {\em model}, $P$ is a distribution on sample space $\cZ$, and
$\loss: \cF \times \cZ \rightarrow \reals \cup \{\infty \}$ is a loss
function.  We denote by $\loss_f(z) \coloneqq \loss(f,z)$ the loss of
predictor $f \in \cF$ under outcome $z \in \cZ$. If $Z \sim P$, we
abbreviate $\loss_f(Z)$ to $\loss_f$. In all our examples,
$\cZ = \cX \times \cY$. We obtain
e.g.\ standard (random-design) regression with squared loss by taking $\cY = \reals$ and $\cF$ to be  some subset of the class of all functions
$f: \cX \rightarrow \reals$
and, for $z=(x,y)$, $\loss_f(x,y) = (y -f(x))^2$; logistic regression is obtained by taking $\cF$ as before, $\cY = \{-1,1\}$ and $\loss_f(x,y) = \log (1 + \exp(- y f(x))$. 
We get conditional density estimation by taking $\{p_f(Y \mid X) :f \in \cF\}$ to be a family of conditional probability mass or density functions (defined relative to some measure $\mu$), extended to $n$ outcomes by the i.i.d.\ assumption, and taking conditional log-loss $\ell_f(x,y) \coloneqq - \log p_f(y \mid x)$.

We are given an i.i.d.\ sample
$Z^n \coloneqq Z_1, Z_2, \ldots, Z_n \sim P$ where each $Z_i$ takes
values in $\cZ$, and we consider, as our learning algorithm, the {\em
  generalized Bayesian posterior}, also known as the {\em Gibbs
  posterior}, $\Pi_n$ on $\cF$, defined by its density
\begin{equation}\label{eq:bayesgenpost}
\pi_n (f) 
\coloneqq  \frac{\exp\left(-\eta \sum_{i=1 }^{n} \loss_f(z_i) \right) \cdot \pi_0(f)}{\int_\model
	\exp\left(-\eta \sum_{i= 1 }^{n} \loss_f(z_i)\right) \cdot \pi_0(f) \text{d} \rho(f)},
\end{equation}
where $\eta>0$ is the \emph{learning rate}, and $\pi_0$ is the density
of some prior distribution $\Pi_0$ on $\cF$ relative to an underlying
measure $\rho$.  Note that, in the conditional log-loss setting, we
get that
\begin{equation}\label{eq:genbayeslik}
  \pi_n (f) \propto \prod_{i=1}^n (p_f(y_i \mid x_i))^{\eta} \pi_0(f),
  \end{equation}
  which, if $\eta=1$, reduces to standard Bayesian inference.
  While GM's result (quoted as Theorem~\ref{thm:metric} below) works for arbitrary loss functions, Theorem~\ref{thm:expfam-exptails} and our empirical simulations (this paper's new results) revolve around (generalized) linear models. For these models, \eqref{eq:bayesgenpost} can be equivalently interpreted either in terms of the original loss functions $\ell_f$ or in terms of the conditional likelihood $p_f$.
  For example, consider regression with $\loss_f(x,y) = (y -f(x))^2$
  and fixed $\eta$. Then (\ref{eq:bayesgenpost}) induces the same
  posterior distribution $\pi_n(f)$ over $\cF$ as does
  (\ref{eq:genbayeslik}) with the conditional distributions
  $p_f(y|x) \propto \exp(-(y-f(x))^2$, which is again the same as
  (\ref{eq:bayesgenpost}) with $\ell_f$ replaced by the conditional
  log-loss $\ell' _f(x,y) \coloneqq - \log p_f(y|x)$, giving a
  likelihood corresponding to Gaussian errors with a particular fixed
  variance; an analogous statement holds for logistic
  regression. Thus, all our examples can be interpreted in terms of
  (\ref{eq:genbayeslik}) for a model that is misspecified, i.e., the
  density of $P(Y|X)$ is not equal to $p_f$ for any $f \in \cF$. As is
  customary (see e.g.~\citet{bartlett2005local}), we assume throughout
  that there exists an optimal $\fopt \in \cF$ that achieves the
  smallest \emph{risk} (expected loss)
  $\E[\loss_{\fopt}(Z)] = \inf_{f \in \cF} \E[\loss_f(Z)]$. If $\cF$
  is a GLM, the risk minimizer again has additional interpretations:
  first, $f^*$ minimizes, among all $f \in \cF$, the conditional KL
  divergence
  ${\bf E}_{(X,Y) \sim P} [\log \left( p(Y|X)/ p_f(Y|X) \right)]$ to
  the true distribution $P$. Second, if there is an $f \in \cF$ with
  ${\bf E}_{X,Y \sim P}[Y \mid X] = f(X)$ (i.e.\ $\cF$ contains the
  {\em true regression function}, or equivalently, {\em true
    conditional mean}), then the risk minimizer satisfies $f^* = f$.
\subsection{Bad Behavior of Standard Bayes}
\label{sec:bad}
\begin{figure}
	\centerline{\includegraphics[width=0.5\textwidth]{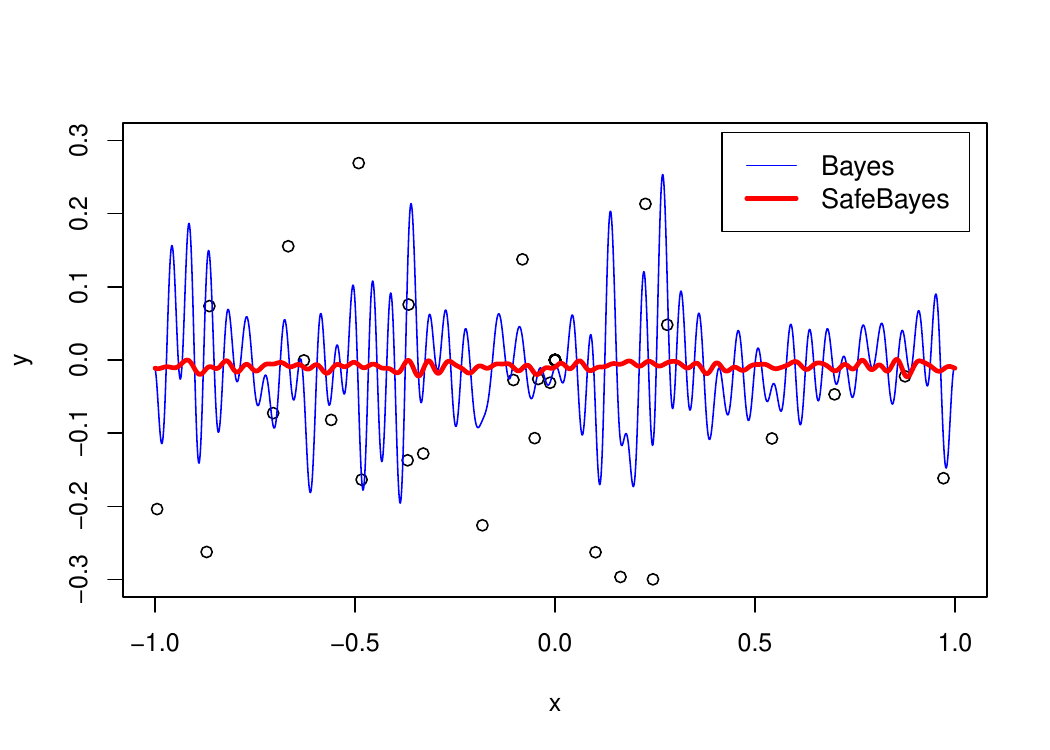}}
	\caption{\emph{Predictions of standard Bayes (blue) and SafeBayes (red), $n=50$, $p=101$.} \vspace*{-0.4 cm} \label{fig:ex1}}
\end{figure}

\begin{example}\label{ex:basic}{\rm 
    We consider a Bayesian lasso regression setting \citep{Park} with
    random design, with a Fourier basis. We sample data $Z_i = (X_i,
    Y_i)$ i.i.d.~$\sim P$, where $P$ is defined as follows: we first sample {\em preliminary\/}
    $(X'_i,Y'_i)$ with $X'_i \overset{i.i.d.}{\sim}$
    Uniform$([-1,1])$; the dependent variable $Y'_i$ is set to $Y'_i = 0
    + \epsilon_i$, with $\epsilon_i \sim \cN(0, \sigma^2)$ for some fixed
    value of $\sigma$, independently of $X'_i$. In other words: the
    true distribution for $(X'_i, Y'_i)$ is `zero with Gaussian
    noise'. Now we toss a fair coin for each $i$. If the coin lands
    heads, we set the actual $(X_i,Y_i) \coloneqq (X'_i,Y'_i)$, i.e.\ we keep
    the $(X'_i, Y'_i)$ as they are, and if the coin lands tails, we
    put the pair to zero: $(X_i, Y_i) \coloneqq (0,0)$. 

We model the relationship between $X$ and $Y$ with a $p^{\mathrm{th}}$ order Fourier basis. Thus, $\cF = \{ f_{\beta}: \beta \in \reals^{2p+1} \}$, with
$f_{\beta}(x)$ given by
\begin{align*}
 \left\langle \beta,  \frac{1}{\pi} \cdot  \left(2^{-1/2}, \cos(x), \sin(x), \cos(2x),  \mydots,  \sin(px)\right) \right\rangle,
\end{align*}
    and the $\eta$-posterior is defined by (\ref{eq:bayesgenpost})
    with $\ell_{f_{\beta}}(x,y) = (y- f_{\beta}(x))^2$; the prior is
    the Bayesian lasso prior whose definition we recall in
    Section~\ref{sec:lassosamp}. Since our `true' regression function
    ${\bf E}[Y_i \mid X_i]$ is $0$, in an actual sample around $50\%$
    of points will be noiseless, \emph{easy} points, lying on the true
    regression function. Since the actual sample of $(X_i,Y_i)$ has less
    noise then the original sample $(X'_i, Y'_i)$, we would 
    expect Bayesian lasso regression to learn the correct
    regression function, but as we see in the blue line in
    Figure~\ref{fig:ex1}, it overfits and learns the noise instead (later on (Figure~\ref{fig:risk-simulaties resultaten} in Section~\ref{sec:lassoandhs}) we shall see that, not surprisingly, this results in terrible predictive behavior). By
    removing the noise in half the data points, we misspecified the
    model: we made the noise heteroscedastic, whereas the model
    assumes homoscedastic noise. Thus, in this experiment the {\em
      model is wrong}. Still, the distribution in $\cF$ closest to the
    true $P$, both in KL divergence and in terms of minimizing the
    squared error risk, is given by the conditional distribution
    corresponding to $Y_i = 0 + \epsilon_i$, where $\epsilon_i$ is
    i.i.d.~$\sim \cN(0,\sigma^2)$. While this element of $\cF$ is in
    fact {\em favored\/} by the prior (the lasso prior prefers $\beta$
    with small $\|\beta\|_1$), nevertheless, for small samples, the
    standard Bayesian posterior puts most if its mass at $f$ with many
    nonzero coefficients.  In contrast, the generalized posterior
    \eqref{eq:bayesgenpost} with $\eta = 0.25$ gives excellent results
    here. To learn this $\eta$ from the data, we can use the
    Safe-Bayesian algorithm of \citet{Grunwald12}. The result is
    depicted as the red line in Figure~\ref{fig:ex1}. Implementation
    details are in Section \ref{sec:lassosamp} and
    Appendix~\ref{app:impl}; the details of the figure are in
    Appendix~\ref{ap:detailsfigures}. }
\end{example}

The example is similar to that of \citet{grunwald2017inconsistency},
who use multidimensional $X$ and a ridge (normal) prior on
$\|\beta\|$; in their  example, standard Bayes succeeds when
equipped with a lasso prior; by using a trigonometric basis we can make
it `fail' after all.  \citet{grunwald2017inconsistency} relate the
potential for the overfitting-type of behavior of standard Bayes, as
well as the potential for full inconsistency (i.e.\ even holding as
$n \rightarrow \infty$) as noted by \citet{GrunwaldL07}
to properties of the Bayesian predictive distribution
${\bar{p}(Y _{n+1} \mid X_{n+1}, Z^n)} \coloneqq {\int_\model p_f(Y_{n+1}
  \mid X_{n+1}) \pi_n(f \mid Z^n) \text{d} \rho(f)}$.  Being a mixture
of $f \in \cF$, $\bar{p}(Y_{n+1} \mid X_{n+1})$, is a member of the
convex hull of densities $\cF$ but not necessarily of $\cF$ itself.
As explained by Gr\"unwald and Van Ommen, severe
overfitting may take place if $\bar{p}(Y_{n+1} \mid X_{n+1}, Z^n)$ is
`far' from any of the distributions in $\cF$. It turns out that this
is exactly what happens in the lasso example above, as we see from
Figure \ref{fig:predvar} (details in
Appendix~\ref{ap:detailsfigures}). 
This figure plots the data points as $(X_i, 0)$ to indicate their location; we see that the predictive variance of standard Bayes fluctuates, being small around the data points and large elsewhere. However, it is obvious that for every density $p_f$ in our model $\cF$, the variance is fixed independently of $X$, 
and thus $\bar{p}(Y_{n+1} \mid X_{n+1}, Z^n)$ is indeed very far from any particular $p_f$ with $f \in \cF$. 
In contrast, for the generalized Bayesian lasso with $\eta = 0.25$, the
corresponding predictive variance is almost constant; thus, at the
level $\eta = 0.25$ the predictive distribution is almost `in-model' (in
machine learning terminology, we may say that $\bar{p}$ is `proper'
\citep{shalev2014understanding}, and the overfitting behavior then does
not occur anymore.
\begin{figure}
	\centerline{\includegraphics[width=0.5\textwidth]{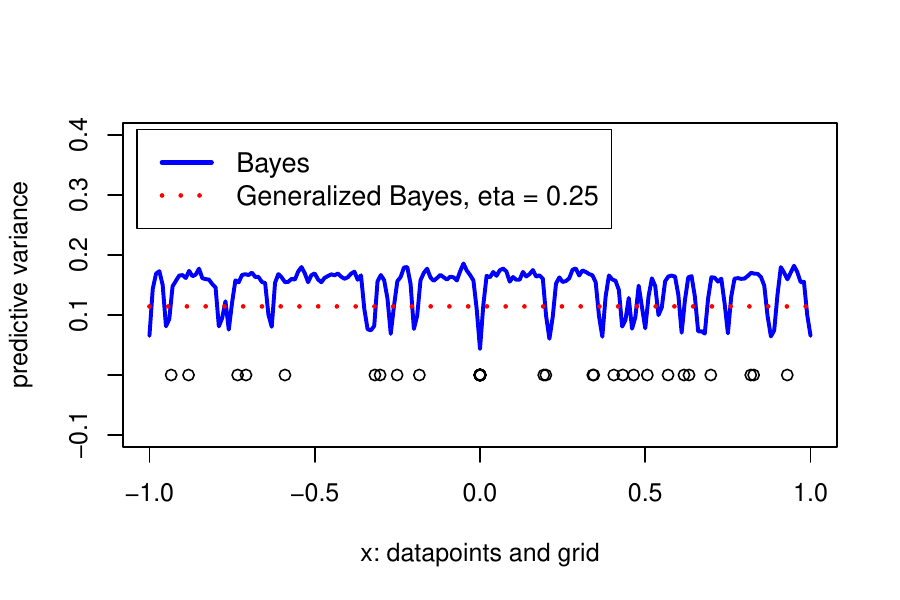}}
	\caption{\emph{Variance of Predictive Distribution $\bar{p}(Y_{n+1} \mid X_{n+1}, Z^n)$ for a single run with $n=50$}. \vspace*{-0.4 cm} \label{fig:predvar}}
\end{figure}
\subsection{When Generalized Bayes Concentrates}
\label{sec:consistency}
Having just seen bad behavior for $\eta=1$, we now recall some
results from GM. Under some conditions, GM show that generalized Bayes, 
for appropriately chosen $\eta$, does concentrate at fast rates
even under misspecification.   
We first recall (a very special case of) the
asymptotic behavior under misspecification theorem of GM.
GM bound (a) the {\em misspecification metric\/} $\dhel_{\bar\eta}$ in
terms of (b) the {\em information complexity}.  The bound (c) holds under a simple condition on the
learning problem that was termed the {\em central condition\/} by
\citet{erven2015fast}. Before presenting the theorem we explain
(a)--(c).
As to (a), we define the {\em
  misspecification metric\/} $\dhel_{\bar\eta}$ in terms of its square
by
	\begin{align*}
\dhel^2_{\bar\eta}(f,f') \coloneqq \frac{2}{\bar{\eta}} \left(1- \int \sqrt{p_{f,\bar\eta}(z) p_{f',\bar\eta}(z)} \text{d} \mu(z) \right)
\end{align*}
which is the ($2 /\bar\eta$-scaled) squared Hellinger distance between
$p_{f,\bar\eta}$ and $p_{f',\bar\eta}$.  Here, a density
$p_{f,\bar{\eta}}$ is defined as
\begin{equation} \nonumber
p_{f,\bar\eta}(z) 
\coloneqq p(z) \frac{\exp(-\bar\eta \xsloss{f}(z))}{\E[\exp(-\bar\eta \xsloss{f}(Z))]},
\end{equation}
where $\xsloss{f} = \loss_f - \loss_{\fopt}$ is the {\em excess
  loss\/} of $f$.
GM show that $\dhel_{\bar\eta}$ defines a metric for all $\bar\eta > 0$.
If $\bar\eta=1$, $\ell$ is log-loss, and the model is well-specified,
then it is straightforward to verify that $p_{f,\bar\eta} = p_f$, 
and so $(1/2) \cdot \dhel_{\bar\eta}$ becomes the standard squared Hellinger distance.

As to (b), we denote by  $\rsc_{n, \eta}(\Pi_0)$ the information complexity, defined as:
\begin{align}
\rsc_{n,\eta}(\Pi_0) 
\coloneqq \Exp_{\rv{f} \sim \dol_n} \left[ \frac{1}{n} \sum_{i=1}^n \xslossat{\rv{f}}{Z_i} \right]
+ \frac{\KL( \dol_n \pipes \Prior)}{\eta \cdot n  } = 
\nonumber \hspace*{-0.5 cm} \\ \  \hspace*{-2 cm} \ \label{eq:bayesmarginal} 
- \frac{1}{ \eta n} \log \int_\model \pi_0(f) e^{- \eta \sum_{i=1}^n \ell_f(Z_i)} \text{d} \rho(f) -
\sum_{i=1}^n \ell_{f^* }(Z_i),
\end{align}
where $\rv{f}$ denotes the predictor sampled from the posterior
$\Pi_n$ and $\KL$ denotes KL divergence; we suppress dependency of $\rsc$ on $\fopt$ in the
notation. The fact that both lines above are equal (noticed by, among
others, \citet{zhang2006information}; GM give an explicit proof)
allows us to write the information complexity in terms of a
generalized Bayesian predictive density which is also known as {\em
  extended stochastic complexity\/} \citep{Yamanishi98}. It also plays
a central role in the field of prediction with expert advice as the
{\em mix-loss\/} \citep{erven2015fast,CesaBianchiL06} and coincides with the minus log
of the standard Bayesian predictive density if $\eta=1$ and $\ell$ is
log-loss. It can be thought of as a complexity measure analogous to
VC dimension and Rademacher complexity.

As to (c), GM's result holds under the \emph{central condition}
(\citep{li1999estimation}; name due to
\cite{erven2015fast}) which expresses that, for some fixed
$\bar\eta > 0$, for all fixed $f$, the probability that the loss of $f$
exceeds that of the optimal $f^*$ by $a/\bar\eta$ is exponentially
small in $a$:
\begin{definition}[Central Condition, Def.~7 of GM]\label{def:central}
	Let $\bar{\eta} > 0$. We say that $\smtuple$ satisfies the 
	\emph{$\bar\eta$-strong central condition} if,
        for all $f \in \cF$: $\E\left[e^{-\bar\eta \xsloss{f}}  \right] \leq 1$.
      \end{definition}
      As straightforward rewriting shows, this condition holds {\em
        automatically\/}, for any $\bar\eta \leq 1$ in the density
      estimation setting, if the model is correct;
      \citet{erven2015fast} provide some other cases in which it holds, and show that many other conditions on $\ell$ and $P$ that allow fast rate convergence that have been considered before in the statistical and on-line learning literature, such as {\em exp-concavity\/} \citep{CesaBianchiL06}, the {\em Tsybakov\/} and {\em Bernstein\/} conditions \citep{bartlett2005local,Tsybakov04} and several others, can be viewed as special cases of the central condition; yet 
they don't discuss GLMs. Here is GM's result:
\begin{theorem}[Theorem~10 from GM]\label{thm:metric} 
	Suppose that the $\bar{\eta}$-strong central condition holds. Then for any $0 < \eta < \bar\eta$, the metric $\dhel_{\bar\eta}$ satisfies
	\begin{align*}
	\E_{Z^n \sim P} \E_{\rv{f} \sim \dol_n} \left[ \dhel^2_{\bar\eta}(\fopt, \rv{f}) \right] 
	\leq  C_{\eta} \cdot \E_{Z^n \sim P}\left[ \rsc_{n,\eta}(\Pi_0) \right]
	\end{align*}
	with $C_{\eta} = \eta / (\bar\eta- \eta)$.
	In particular, $C_{\eta} < \infty$ for $0 < \eta < \bar\eta$, and $C_{\eta} = 1$ for $\eta = \bar\eta/2$.
      \end{theorem}
      Thus, we expect the posterior to concentrate at a rate dictated
      by ${\bf E}[\rsc_{n,\eta}]$ in neighborhoods of the best
    (risk-minimizing, KL optimal, or even true regression function)
    $f^*$.  
    The misspecification metric
    $\dhel^2_{\bar\eta}$ on the left hand side is a weak metric,
    however, in Appendix~\ref{sec:witness} we show that we can replace it by
    stronger notions such as KL-divergence, squared error or logistic
    loss.  
    Theorem \ref{thm:metric} generalizes previous results (e.g.~\citet{zhang2006epsilon, zhang2006information}) to the misspecified setting. In the well-specified case,
    Zhang, as well as several other authors
    \citep{WalkerH02,martin2017empirical}, state a result that holds for
    any $\eta < 1$ but not $\eta=1$. This suggests that there is an advantage 
    to taking $\eta$ slightly smaller than one even when the model is
    well-specified (for more details see \citet{zhang2006epsilon}).

      To make the theorem work for GLMs under misspecification, we
      must verify (a) that the central condition still holds (which is
      in general not guaranteed) and that (b) the information
      complexity is sufficiently small.  As to (a), in the following
      section we show that the central condition holds (with
      $\bar\eta$ usually $\neq 1$) for $1$-dimensional exponential
      families and high-dimensional generalized linear models (GLMs)
      if the noise is misspecified, as long as $P$ has exponentially
      small tails; in particular, we relate $\bar\eta$ to the variance of $P$. 
      As to (b), if the model is correct (the conditional distribution
      $P(Y \mid X)$ has density $f$ equal to $p_f$ with $f \in \cF$),
      where $\cF$ represents a $d$-dimensional GLM, then it is known (see e.g.~\citet{zhang2006information}) that, for any
      prior $\Pi_0$ with continuous, strictly positive density on
      $\cF$, the information complexity satisfies
\begin{equation}\label{eq:bic}
   {\bf E}_{Z^n \sim P} \left[   \rsc_{n,\eta}(\Pi_0)\right] = O \left( \frac{d}{n} \cdot \log n
      \right),
\end{equation}
      which leads to bounds within a log-factor of the minimax optimal
      rate (among all possible estimators, Bayesian or not), which is $O(d/n)$. 
      While such results were only
      known for the well-specified case, in
      Proposition~\ref{prop:tenerife} below we show that, for GLMs, they
      continue to hold for the misspecified case.
     
\section{GENERALIZED GLM BAYES}
\label{sec:genglm}
Below we first show that the central condition holds for natural univariate
exponential families; we then extend this result to the GLM case, and
establish bounds in information complexity of GLMs.  Let the class
$\cF = \{ p_{\theta}: \theta \in \Theta\}$ be a univariate natural
exponential family of distributions on $\cZ = \cY$, represented by
their densities, indexed by natural parameter $\theta \in \Theta
\subset \reals$ \citep{BarndorffNielsen78}.  The elements of this
restricted family have probability density functions
\begin{align}\label{eq:expfamdef}
p_\theta(y) \coloneqq \exp( \theta y - F(\theta) + r(y) ) ,
\end{align}
for log-normalizer $F$ and carrier measure $r$. We denote the
corresponding distribution as $P_{\theta}$. In the first part of the theorem below we assume that $\Theta$ is restricted to an arbitrary closed
interval $[\underline\theta,\bar\theta]$ with
$\underline\theta < \bar\theta$ that resides in the interior of the
natural parameter space
$\bar{\Theta} = \{ \theta: F(\theta) < \infty\}$. Such $\Theta$
allow for a simplified analysis because within $\Theta$ the
log-normalizer $F$ as well as all its derivatives are uniformly
bounded from above and below; see \eqref{eq:ineccsi} in
Appendix~\ref{app:proofs}.  As is well-known (see
e.g.~\citet{BarndorffNielsen78}), exponential families can
equivalently be parameterized in terms of the mean-value
parameterization: there exists a $1$-to-$1$ strictly increasing
function $\mu: \bar{\Theta} \rightarrow \reals$ such that $\E_{Y \sim
  P_{\theta}}[Y] = \mu(\theta)$.  As is also well-known, the density
$p_{\fopt} \equiv p_{\theta^*}$ within $\cF$ minimizing KL divergence
to the true distribution $P$ satisfies $\mu(\theta^*) = {\bf E}_{Y\sim
  P}[Y]$, whenever the latter quantity is contained in $\mu(\Theta)$
\citep{Grunwald07}.  In words, the best approximation to $P$ in $\cF$
in terms of KL divergence has the same mean of $Y$ as $P$.

\begin{theorem} \label{thm:expfam-exptails}
  Consider a learning problem $\smtuple$ with $\loss_{\theta}(y) = -
  \log p_{\theta}(y)$ the log loss and $\cF= \{p_{\theta}: \theta \in
  \Theta\}$ a univariate exponential family as above. \\ (1). Suppose
  that $\Theta = [\underline\theta,\bar\theta]$ is compact as above
  and that $\theta^* = \arg\min_{\theta \in \bar\Theta} D(P\|
  P_{\theta})$ lies in $\Theta$. Let $\sigma^2> 0$ be the true
  variance $\E_{Y \sim P}(Y - \E[Y])^2$ and let $(\sigma^*)^2$ be the
  variance $\E_{Y \sim P_{\theta^*}}(Y - \E[Y])^2$ according to
  $\theta^*$.  Then \subitem(i) for all $\bar\eta >
  (\sigma^*)^2/\sigma^2$, the $\bar\eta$-central condition does {\em
    not\/} hold. 
 \subitem(ii) Suppose there exists $\eta^{\circ} > 0$
  such that {$\bar{C} \coloneqq \E_P[\exp(\eta^{\circ} |Y|)] < \infty$.}
  Then there exists $\bar\eta > 0$, depending only on $\eta^{\circ}$,
  $\bar{C}, \underline{\theta}$ and $\overline{\theta}$
  such that the $\bar\eta$-central condition holds. Moreover, \subitem(iii),
  for all $\delta > 0$, there is an $\epsilon > 0$ such that, for all
  $\bar\eta \leq (\sigma^*)^2/\sigma^2 - \delta$, the
  $\bar{\eta}$-central condition holds relative to the restricted
  model $\cF_{\epsilon} = \{p_{\theta}: \theta \in [\theta^*-
    \epsilon,\theta^*+\epsilon]$\}\vspace*{-0.2 cm}.
        \end{theorem}
{\em (2).
          Suppose that $P$ is Gaussian with variance $\sigma^2> 0$ and that
		$\cF$ indexes a full Gaussian location family. Then
		the $\bar\eta$-central condition holds iff $\bar\eta \leq
		(\sigma^*)^2/\sigma^2$.
              }
              
We provide (iii) just to give insight --- `locally', i.e.~in
restricted models that are small neighborhoods around the
best-approximating $\theta^*$, the smallest $\bar\eta$ for which
the central condition holds is determined by a ratio of variances. The final part
shows that for the Gaussian family, the same holds not just locally
but globally (note that we do not make the compactness assumption on
$\Theta$ there); we warn the reader though that the standard posterior ($\eta=1$) based on a model with fixed variance $\sigma^*$ is quite different from the generalized posterior with $\eta = (\sigma^*)^2/\sigma^2$ and a model with variance $\sigma^2$ \citep{grunwald2017inconsistency}. Finally, while in practical cases we often find $\bar\eta < 1$ (suggesting that Bayes may only succeed if we learn `slower' than with the standard $\eta =1$, i.e.\ the prior becomes more important), the result shows that we can  also very well have $\bar\eta > 1$; we give a practical example at the end of Section~\ref{sec:examples}.
\commentout{
This
happens if {\em data from $P$ provide more information about
  $\theta^*$ than would be expected if $P_{\theta^*}$ itself were
  true}.  As an extreme case, suppose, for example, that our model is
Poisson, geometric, or any other exponential family supported on
$\naturals$, and the true data generating process is deterministic,
always producing the same number - so our data might be,
e.g.~$Y^n = 4,4,\ldots,4$. Then the larger we pick $\bar\eta$, the
faster the $\bar\eta$-posterior concentrates around the optimal value
$\theta^*$ with $\mu(\theta^*) = 4$.
}
Theorem~\ref{thm:expfam-exptails} is new and supplements \citeauthor{erven2015fast}'s (\citeyear{erven2015fast}) various examples of $\cF$ which
satisfy the central condition.  In the
theorem we require that both tails of $Y$ have exponentially small
probability.

\paragraph{Central Condition: GLMs}

Let $\cF$ be the generalized linear model \citep{McCullaghN89} (GLM) indexed by parameter $\beta \in \cB
\subset \reals^d$ with link function $g: \reals \rightarrow
\reals$. By definition this means that there exists a set $\cX \subset
\reals^{d}$ and a
univariate exponential family $\cQ = \{ p_{\theta} : \theta \in
\bar\Theta \}$ on $\cY$ of the form \eqref{eq:expfamdef} such that the
conditional distribution of $Y$ given $X=x$ is, for all possible values
of $x \in \cX$, a member of the family $\cQ$, with mean-value parameter
$g^{-1}(\langle \beta, x  \rangle)$.
Then
the class $\cF$ can be written as $\cF = \{ p_\beta: \beta \in
\mathcal{B}\}$, a set of conditional probability density functions
such that
\begin{align}
p_\beta(y \mid x) 
\coloneqq \exp \bigl( \theta_x(\beta) y - F(\theta_x(\beta)) + r(y) \bigr),
\label{eq:toulon}
\end{align}
where $\theta_x(\beta) \coloneqq \mu^{-1}(g^{-1}(\langle \beta, x  \rangle))$,
and $\mu^{-1}$, the inverse of $\mu$ defined above,  sends mean parameters to
natural parameters.  We then have $\E_{P_\beta} [ Y \mid X ] =
g^{-1}(\langle \beta, X  \rangle)$, as required.

\begin{proposition}\label{prop:tenerife} Under the following three assumptions, the learning problem $\smtuple$ with $\cF$ as above satisfies the $\bar\eta$-central condition for
some $\bar\eta>0$ depending only on the parameters of the problem:
\begin{enumerate}
	\item (Conditions on $g$): 
	the inverse link function $g^{-1}$ has bounded derivative on the domain $\cB \times \cX $, and the image of the inverse link on the same domain is a bounded interval in the interior of the mean-value
	parameter space
	$\{ \mu\in \reals: \mu = \E_{Y \sim q}[Y]\; :\; q \in \cQ \}$ (for all standard link functions, this can be enforced by restricting $\cB$ and $\cX$ to an (arbitrarily large but still) compact domain).
	\item (Condition on `true' $P$): for some $\eta > 0$ we have \\
	${\sup}_{x \in \cX} \E_{Y \sim P}[\exp(\eta | Y|) \mid X=x] < \infty$.
	\item (Well-specification of  conditional mean): there exists $\beta^{\circ} \in \cB$ such that $\E [ Y \mid
	X ] = g^{-1}(\langle \beta^{\circ}, X  \rangle)$.
\end{enumerate}
\end{proposition}
A simple argument (differentiation with respect to $\beta$) shows that
under the third condition, it must be the case that
$\beta^{\circ} = \beta^*$, where $\beta^* \in \cB$ is the index
corresponding to the density $p_{\fopt} \equiv p_{\beta^*}$ within
$\cF$ that minimizes KL divergence to the true distribution $P$.
Thus, our conditions imply that $\cF$ contains a $\beta^*$ which
correctly captures the conditional mean (and this will then be the
risk minimizer); thus, as is indeed the case in
Example~\ref{ex:basic}, the regression function must be well-specified
but the noise can be severely misspecified.

We stress that the three conditions have very different statuses. The
first is mathematically convenient; it can be enforced by truncating
parameters and data, which is awkward but may not lead to substantial
deterioration in practice.  Whether it is even really needed or not
is not clear (and may in fact depend on the chosen exponential
family).  The second condition is really necessary --- as can
immediately be seen from Definition~\ref{def:central}, the strong
central condition cannot hold if $Y$ has polynomial tails and for some
$f$ and $x$, $\loss_f(x,Y)$ increases polynomially in $Y$ (in
Section~6 of their paper, GM consider weakenings of the central
condition that still work in such situations). For the third
condition, however, we suspect that there are many cases in which it
does not hold yet still the strong central condition holds; so then
the GM convergence result would still be applicable under `full
misspecification'; investigating this will be the subject
of future work.
\paragraph{GLM Information Complexity}
To apply Theorem~\ref{thm:metric} to get convergence bounds for
exponential families and GLMs, we need to verify that the central
condition holds (which we just did) and we need to bound the
information complexity, which we proceed to do now. It turns out that
the bound on $\rsc_{n,\eta}$ of $O( (d/n) \log n)$ of (\ref{eq:bic})
continues to hold unchanged under misspecification, as is an immediate
corollary of applying the following proposition to the definition of
$\rsc_{n,\eta}$ given above (\ref{eq:bayesmarginal}):
\begin{proposition}\label{prop:entroboundb}
	Let $\smtuple$ be a learning problem with $\cF$ a GLM satisfying
	Conditions 1--3 above. Then for all $f \in \cF$, 
	$ \E_{X,Y \sim P}[L_f] = \E_{X,Y \sim P_{\fopt}} [L_f].
	$
\end{proposition}
This result follows almost immediately from the `robustness property
of exponential families' (Chapter 19 of \citet{Grunwald07}); for
convenience we provide a proof in Appendix~\ref{app:proofs}.
The result implies that any bound in $\rsc_{n,\eta}(\Pi_0)$ for a
particular prior in the well-specified GLM case, in particular
(\ref{eq:bic}), immediately transfers to the same bound for the
misspecified case, as long as our regularity conditions hold, allowing
us to apply Theorem~\ref{thm:metric} to obtain the parametric rate
for GLMs under misspecification.  \commentout{ Here is a concrete
  example combining Proposition~\ref{prop:entroboundb} and
  \eqref{eq:ggvc} (see also Eq.~(25) of GM).
Consider a finite $(M < \infty)$ or countably infinite $(M = \infty)$
union $\cF \coloneqq \bigcup_{j= 1, \ldots, M} \cF_j$ of nested GLMs, i.e.~$\cF_1
\subset \cF_2 \subset \cF_3 \subset \ldots$ with common domain $\cX$
and link function $g$. Here $\cF_j$ has parameter space $\cB_j \subset
\reals^{d_j}$ for a strictly increasing sequence of dimensions $d_1,
d_2, \ldots$ such that each $\cF_j$ satisfies Condition 1 underneath \eqref{eq:toulon}. For example, we could have that there is an underlying 1-dimensional covariate $U$, and  $X = (f_1(U), f_2(U), f_d(U))$ are the first $d$ components of some basis expansion of $U$ such as $f_d = U^d$. We equip each $\cF_j$ with an arbitrary prior density
$\pi^{(j)}$ bounded away from $0$, which is possible because of the
imposed Condition 1. We construct a prior on $\cF$ by taking an arbitrary
discrete prior $\pi'$ with full support on $\{1, \ldots, M\}$, setting $\pi(f,
\cF_j) \coloneqq \pi^{(j)}(f) \pi'(j)$.

Now suppose that $\smtuple$ is such that $\cF$ satisfies Condition
1--3 underneath \eqref{eq:toulon}. Condition 3 implies existence of an
$\fopt \in \cF$ that achieves $\min_{f \in \cF} \E[\loss_f]$. Let $k^*$
be the smallest $k$ such that $\cF_k$ contains $\fopt$. Thus, for
simplicity, we only consider the pseudo-nonparametric case where the
model may be nonparametric but the truth is contained in an (unknown)
parametric submodel $\cF_{k^*}$. Proposition~\ref{prop:entroboundb}
implies that in the misspecified case we can still get the
well-specified parametric rate:
\begin{proposition}\label{prop:tenerife}
	Under the conditions on $\cF$ imposed above, the
	$\bar\eta$-central condition holds for some $\bar\eta > 0$. For
	every fixed $\eta> 0$, in particular for $0 < \eta  < \bar\eta$, the $\eta$-generalized Bayesian estimator
	$\Pi_n$ based on a prior of the form specified above satisfies
	$\E[\rsc_{n,\eta}(\Pi_0)] = \tilde{O}(d_{k^*}/n)$, where $k^*$ is the dimension of the submodel 
\end{proposition}
Presumably, this result can be extended to obtain rates for the
nonparametric case where the underlying function is not in the union
of the $\cF_k$ but satisfies appropriate smoothness conditions,
analogously to \citet{barron1991minimum}.
}
\section{MCMC SAMPLING}\label{sec:sampling}

Below we devise MCMC algorithms for obtaining samples from the $\eta$-generalized posterior distribution for two problems: regression and classification. In the regression context we consider one of the most commonly used sparse parameter estimation techniques, the lasso. For classification we use the logistic regression model. In our experiments in Section~\ref{sec:examples}, we compare the performance of generalized Bayesian lasso with Horseshoe regression \citep{horseshoe}. The derivations of samplers are given in Appendix~\ref{app:impl}.

\subsection{Bayesian lasso regression}
\label{sec:lassosamp}
Consider the regression model $Y= X \beta+\eps$,
where $\beta\in\RR^p$ is the vector of parameters of interest, \mbox{$Y\in\RR^n$}, \mbox{$X\in\RR^{n\times p}$}, and \mbox{$\eps\sim \cN(0,\sigma^2 I_n)$} is a noise vector.
The Least Absolute Shrinkage and Selection Operator (LASSO) of \citet{Lasso} is a regularization method used in regression problems for shrinkage and selection of features. The lasso estimator is defined as
$
{\hat{\beta}}_{\text{lasso}} \coloneqq \argmin_{{\beta}} \| Y -  X{\beta}\|_2^2 + \lambda \|\beta\|_1\, ,
$
where $\|\cdot\|_1,\|\cdot\|_2$ are $l_1$ and $l_2$ norms correspondingly.
It can be interpreted as a Bayesian posterior mode (MAP) estimate when the priors on $\beta$ are given by  independent Laplace distributions. As discovered by \citet{Park}, the same posterior on $\beta$ is also obtained by the following Gibbs sampling scheme: set $\eta = 1$ and
denote $ {D_\tau}\coloneqq\diag(\tau_1,\dots,\tau_n).$ Also, let $a\coloneqq\frac{\eta}{2}(n-1) + \frac p2 + \alpha$ and $b_\tau\coloneqq\frac{\eta}{2}( {{Y}-X\beta})^T  ( {{Y}-X\beta}) + \frac12 {\beta}^T  {D_\tau}^{-1}  {\beta} + \gamma$, where $\alpha,\gamma>0$ are hyperparameters. Then the Gibbs sampler is constructed as follows. 
\begin{align*}
\begin{split}
 {\beta} \sim &\,\cN \left( \eta M_\tau {X}^T {{Y}}, \sigma^2M_\tau \right), 
\end{split} \\ 
\begin{split}
\sigma^2 \sim &\,\text{Inv-Gamma} \left( a, b_\tau \right),
\end{split} \quad
{\tau^{-2}_j} \sim \,\text{IG}\left(\displaystyle \sqrt{{\lambda^2\sigma^2}/{\beta_j^2}}, \lambda^2\right),
\end{align*}
where IG is the inverse Gaussian distribution and
$M_\tau\coloneqq(\eta {X}^T {X} + {D_\tau}^{-1})^{-1}$. Following
\citet{Park}, we put a Gamma prior on the shrinkage parameter
$\lambda$. Now, in their paper \citeauthor{Park} only give the scheme for
$\eta = 1$, but, as is straightforward to derive from their paper, the
scheme above actually gives the $\eta$-{\em generalized\/} posterior
corresponding to the lasso prior for general $\eta$ (more details in
Appendix~\ref{app:impl}).  We will use the Safe-Bayesian algorithm for
choosing the optimal $\eta$ developed by
\citet{grunwald2017inconsistency} (see Appendix~\ref{ap:safebayesimpl}). The code for Generalized- and Safe-Bayesian lasso regression can be found in the CRAN $\texttt{R}$-package `SafeBayes' \citep{heide2016safebayes}.

\paragraph{Horseshoe estimator}
The Horseshoe prior is the state-of-the-art global-local shrinkage prior for tackling high-dimensional regularization, introduced by \citet{horseshoe}. Unlike the Bayesian lasso, it has flat Cauchy-like tails, which allow strong signals to remain unshrunk a posteriori. For completeness we include the horseshoe in our regression comparison, using the implementation of \citet{vanderpas2016horseshoe}.

\subsection{Bayesian logistic regression}

Consider the standard logistic regression model $\{ f_{\beta}: \beta \in \reals^{p}\}$, the data $Y_1,\dots, Y_n\in\{0,1\}$ are independent binary random variables observed at the points $X:=(X_1,\dots,X_n)\in \mathbb{R}^{n\times p}$ with
$
P_{f_{\beta}}(Y_i=1 \mid X_i) \coloneqq p_{f_{\beta}}(1 \mid X_i) \coloneqq \frac{e^{X_i^T\beta}}{1+e^{X_i^T\beta}} \, .
$
 The standard Bayesian approach involves putting a Gaussian prior on the parameter $\beta\sim \cN(b,B)$ with mean $b\in\RR^p$ and the covariance matrix $B\in\RR^{p\times p}$. To sample from the $\eta$-generalized posterior we modify a P{\'o}lya--Gamma latent variable scheme described in \citet{polson2013bayesian}. We first introduce latent variables $\omega_1,\dots,\omega_n\in \RR$, which will be sampled from P{\`o}lya-Gamma distribution (constructed to yield a simple Gibbs sampler for Bayesian logistic regression, for more details see \citet{polson2013bayesian}). Let $\Omega:=\diag\{\omega_1,\dots,\omega_n\}$, $\kappa\coloneqq(Y_1-1/2,\dots, Y_n-1/2)^T$, ${V_\omega\coloneqq(X^T\Omega X+B^{-1})^{-1}}$, and $ {m_\omega\coloneqq V_\omega(\eta X^T\kappa+B^{-1}b)}$.
Then the Gibbs sampler for $\eta$-generalized posterior is given by
$
\omega_i\sim \text{PG}(\eta, X_i^T\beta), \quad
\beta\sim \cN(m_\omega, V_\omega), 
$
where PG is the P{\`o}lya-Gamma distribution. 

\section{EXPERIMENTS}\label{sec:examples}
Below we present the results of experiments that compare the performance of the derived Gibbs samplers with their standard counterparts. More details/experiments are in Appendix~\ref{ap:detailsfigures}.

\subsection{Simulated data}\label{sec:lassoandhs}

%
%
 \paragraph{Regression} In our experiments we focus on prediction, and we run simulations to determine the \emph{square-risk} (expected squared error loss) of our estimate relative to the underlying distribution $P$:
\mbox{$\E_{(X,Y) \sim P} (Y - X\beta )^2$},
where $X\beta$ would be the conditional expectation, and thus the square-risk minimizer, if $\beta$ would be the true parameter (vector).

\begin{figure}
	\centerline{\includegraphics[width=0.45\textwidth, height=0.3\textwidth]{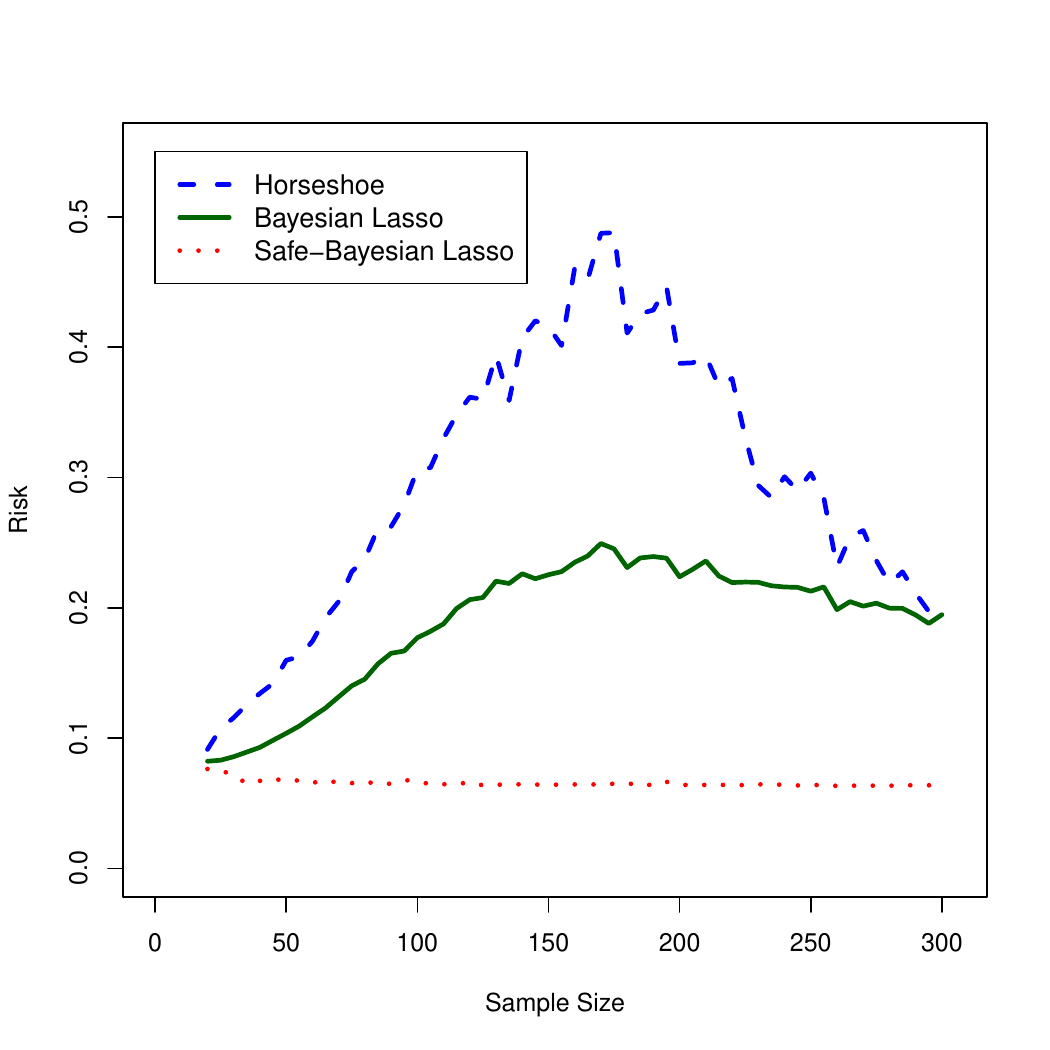}}
	\caption{\emph{Simulated squared error risk (test error) with
          respect to $P$ as function of sample size for the
          \emph{wrong-model} experiments of Section~\ref{sec:lassoandhs} using the posterior predictive
          distribution of the standard Bayesian lasso (green, solid),
          the Safe-Bayesian lasso (red, dotted), both with standard
          improper priors, and the Horseshoe (blue, dashed); and $201$
          Fourier basis functions.} \vspace*{-0.6 cm} \label{fig:risk-simulaties resultaten}}
\end{figure}

Consider the data generated as described in Example~\ref{ex:basic}.
We study the performance of the $\eta$-generalized Bayesian lasso with $\eta$ chosen by the Safe-Bayesian algorithm (we call it the Safe-Bayesian lasso) in comparison with two popular estimation procedures for this context: the Bayesian lasso (which corresponds to $\eta$=1), and the Horseshoe method. In Figure \ref{fig:risk-simulaties resultaten}  the simulated square-risk is plotted as a function of the sample size for all three methods. We average over enough samples so that the graph appears to be smooth ($25$ iterations for SafeBayes, $1000$ for the two standard Bayesian methods). It shows that both the standard Bayesian lasso and the Horseshoe perform significantly worse than the Safe-Bayesian lasso. Moreover we see that the risks for the standard methods initially grows with the sample size (additional experiments not reported here suggest that Bayes will `recover' at very large $n$). 

\begin{figure}
	\centerline{\includegraphics[width=0.4\textwidth, height=0.3\textwidth]{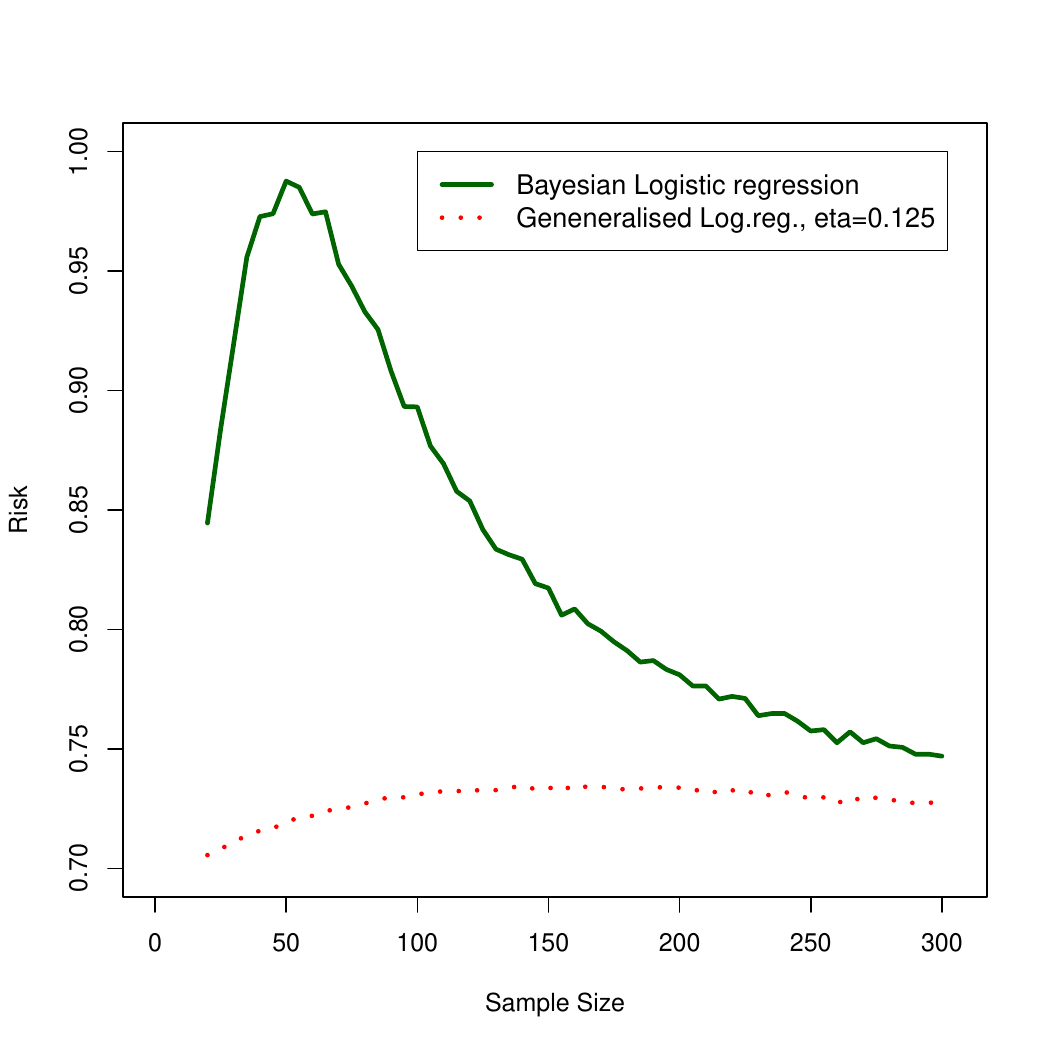}}
	\caption{\emph{Simulated logistic risk as function of  sample size for \emph{wrong-model} experiments of Section \ref{sec:logisticregression} using posterior predictive distribution of standard Bayesian logistic regression (green, solid), and generalized Bayes ($\eta=0.125$, red, dotted) with $25$ noise dimensions.}\vspace*{-0.6 cm} \label{fig:logregrisk}}
\end{figure}

\paragraph{Classification}\label{sec:logisticregression}
We focus on finding coefficients $ \beta$ for prediction, and our error measure is the expected logarithmic loss, which we call \emph{log-risk}: \mbox{$ \E_{(X,Y) \sim P} \left[- \log \text{Li}_\beta(Y\given X) \right]$},
where $ \text{Li}_\beta(Y\given X)\coloneqq {e^{YX^T\beta}}/({1+e^{X^T\beta}})$.
We start with an example that is very similar to the previous one. We generate a $n\times p$ matrix of independent standard normal random variables with $p=25$. For every feature vector $ X_i$ we sample a corresponding $Z_i \sim \cN(0, \sigma^2)$, as before, and we misspecify the model by putting approximately half of the $Z_i$ and the corresponding $X_{i,1}$ to zero. Next, we sample the labels $Y_i \sim \text{Binom}(\exp(Z_i) / (1 + \exp(Z_i))$. We compare standard Bayesian logistic regression ($\eta=1$) to a generalized version ($\eta = 0.125$). In Figure~\ref{fig:logregrisk} we plot the log-risk as a function of the sample size. As in the regression case, the risk for standard Bayesian logistic regression ($\eta=1$) is substantially worse than the one for generalized Bayes ($\eta=0.125$). Even for generalized Bayes, the risk initially goes up a little bit, the reason being that the prior is {\em too good\/}: it is strongly concentrated around the risk-optimal $\beta^* =0$. Thus, the first prediction made by the Bayesian predictive distribution coincides with the optimal $(\beta=0)$ prediction, and in the beginning, due to noise in the data, predictions will first get slightly worse. This is a phenomenon that also applies to standard Bayes with well-specified models; see for example \cite[Example 3.1]{GrunwaldH04}.

Even for the well-specified case it can be beneficial to use $\eta\neq 1$. It is easy to see that the maximum {\it a posteriori} estimate for generalized logistic regression corresponds to the ridge logistic regression method (which penalizes large $\|\beta\|_2$) with the shrinkage parameter $\lambda=\eta^{-1}$. However, when the the prior mean is zero but the risk minimizer $\beta^*$ is far from zero, penalizing large norms of $\beta$ is inefficient, and we find that the best performance is achieved with $\eta>1$. 

\subsection{Real World Data}\label{sec:real}
We present two examples with real world data to demonstrate that bad
behavior under misspecification also occurs in practice. For these
data sets, we compare the performance of Safe-Bayesian lasso and standard
Bayesian lasso.  As the first example we consider the data of the
daily maximum temperatures at Seattle Airport as a function of the
time and date (source: $\texttt{R}$-package $\texttt{weatherData}$,
also available at $\texttt{www.wunderground.com}$). A second example
is London air pollution data (source: $\texttt{R}$-package
$\texttt{Openair}$, for more details see \citet{openair1,
  openair2}). Here the quantity of interest is the concentration of
nitrogen dioxide (NO$_2$), again as a function of time and date. In
both settings we divide the data into a training set and a test set and
focus on the prediction error. In both examples, SafeBayes picks an
$\hat\eta$ strictly smaller than one. Also, for both data sets the
Safe-Bayesian lasso clearly outperforms the standard Bayesian lasso and the
Horseshoe in terms of mean square prediction error, as  seen
from Table~\ref{table:seattleandlondon} (details in
Appendix~\ref{ap:detailsfigures}).

\begin{table}
	\begin{tabularx}{0.48\textwidth}{| l | X | X | X |}
		\hline
		& Horse-shoe & Bayesian lasso & SafeBayes lasso \\
		\hline  
		 MSE ($(^{\circ} \text{C})^2$) & $6.53$ & $6.16$ & $6.04$ \\
		\hline
		MSE ($(\text{ppm})^2$)  & $1169$ & $1201$ & $1142$ \\
		\hline
	\end{tabularx}
	\caption{\emph{Mean square errors for predictions on the Seattle and London data sets of Section~\ref{sec:real}.}\vspace*{-5 mm} \label{table:seattleandlondon}}
\end{table}


\section{FUTURE WORK}
We provided both theoretical and empirical evidence that
$\eta$-generalized Bayes can significantly outperform standard Bayes
for GLMs. However, the empirical examples are only given for Bayesian
lasso linear regression and logistic regression. In future work we
would like to devise generalized posterior samplers for other GLMs
and speed up the sampler for generalized Bayesian logistic
regression, since our current implementation is slow and (unlike our
linear regression implementation) cannot deal with high-dimensional
(and thus, real-world) data yet. Furthermore, the Safe-Bayesian algorithm
of \cite{Grunwald12}, used to learn $\eta$, enjoys good theoretical
performance but is computationally very slow. Since learning $\eta$
for which the central condition holds (preferably the largest possible
value, since small values of $\eta$ mean slower learning) is essential
for using generalized Bayes in practice, there is a necessity for
speeding up SafeBayes or finding an alternative. A potential solution
might be using cross-validation to learn $\eta$, but its theoretical
properties (e.g.\ satisfying the central condition) are yet to be
established.


\subsubsection*{Acknowledgements}

The project leading to this work has received funding from the European Research Council (ERC) under the European Union’s Horizon 2020 research and innovation programme (grant agreement No 834175).
%
%
\DeclareRobustCommand{\VANDER}[3]{#3}
\bibliography{safebayesbib}

\newpage
\onecolumn
\DeclareRobustCommand{\VANDER}[3]{#2}
\appendix
\section{OUTLINE}\label{app:outline}
The appendix of this paper is organized as follows:
\begin{itemize}
	\item Appendix~\ref{app:proofs} provides the proofs for Section~\ref{sec:genglm}.
	\item Appendix~\ref{sec:witness} shows how we can replace  $\dhel^2_{\bar\eta}$ in Theorem~\ref{thm:metric} by stronger notions.
	\item Appendix~\ref{app:learningrate>1} provides an example in which taking a learning rate larger than results in faster learning under misspecification than $\eta = 1$. 
	\item Appendix~\ref{app:impl} provides (implementation) details on the $\eta$-generalized Bayesian lasso and logistic regression; and the Safe-Bayesian algorithm. 
	\item Appendix~\ref{ap:detailsfigures} contains details for the experiments and figures in the main text, and provides additional figures.
\end{itemize}

\section{PROOFS}
\label{app:proofs}
\subsection{Proof of Theorem~\ref{thm:expfam-exptails}} The
	second part of the theorem about the Gaussian location family
        is a straightforward calculation, which we omit. As to the
        first part (Part (i)---(iii)), we will repeatedly use the
        following fact: for every $\Theta$ that is a nonempty compact
        subset of the interior of $\bar{\Theta}$, in particular for
        $\Theta = [\underline{\theta},\bar{\theta}]$ with
        $\underline\theta < \bar\theta$ both in the interior of
        $\overline{\Theta}$, we have:
\begin{equation}\label{eq:ineccsi}
\begin{aligned}
- \infty <  \inf_{\theta \in \Theta} F(\theta) &< \sup_{\theta \in \Theta} F(\theta) < \infty \\
- \infty < \inf_{\theta \in \Theta} F'(\theta) &< \sup_{\theta \in \Theta} F'(\theta) < \infty \\
0 <  \inf_{\theta \in \Theta} F''(\theta) &< \sup_{\theta \in \Theta} F''(\theta) < \infty.
\end{aligned}
\end{equation}
Now, let $\theta, \theta^* \in \Theta$.
	We can write
	\begin{align}\label{eq:basecamp}
	\E \left[ e^{-\eta (\xslosstheta)} \right]  
	= \E_{Y \sim P} \left[ \left( \frac{p_\theta(Y)}{p_{\theta^*}(Y)} \right)^\eta \right]
	 = \exp \left (-G(\eta(\theta- \theta^*))
	+ \eta F(\theta^*) - \eta F(\theta) \right).
	\end{align}
	where $G(\lambda) =  - \log \E_{Y \sim P} \left[ \exp( \lambda Y)\right]$.
	If this quantity is $- \infty$ for all $\eta > 0$, then (i) holds trivially. 
	If not, then (i) is implied by the following statement: 
		\begin{equation}\label{eq:limitvariance}
		\limsup_{\epsilon \rightarrow 0} \left\lbrace \eta: \text{for all $\theta \in [\theta^*- \epsilon,\theta^* + \epsilon]$},  
		\ \E[\exp(\eta \xsloss{p_{\theta}})] \leq 1 \right\rbrace = \frac{(\sigma^*)^2}{\sigma^2}.
		\end{equation}
                Clearly, this statement also implies (iii). To prove
                (i), (ii) and (iii), it is thus sufficient to prove
                (ii) and (\ref{eq:limitvariance}). We prove both by a
                second-order Taylor expansion (around $\theta^*$) of
                the right-hand side of \eqref{eq:basecamp}.
	
	{\em Preliminary Facts}. By our assumption there is a $
        \eta^{\circ}> 0 $ such that $\E[\exp(\eta^\circ |Y|)] =
        \bar{C} < \infty$.  Since $\theta^* \in \Theta =
            [\underline\theta,\overline\theta]$ we must have for every
            $0 < \eta < \eta^{\circ}/(2 |\overline\theta -
            \underline\theta|)$, every $\theta \in \Theta$,
%
	\begin{equation}\label{eq:finite}
	\begin{aligned}
	  \E[ \exp(2 \eta (\theta - \theta^*) \cdot Y)] \leq
          \E[ \exp(2 \eta | \theta - \theta^*| \cdot |Y|)] \leq 
      \E[ \exp(\eta^{\circ} (| \theta - \theta^*|/|\overline\theta -
        \underline\theta|) \cdot |Y|)] \leq \bar{C} < \infty.
	\end{aligned}
	\end{equation}
	%
	%
	The first derivative of the right of \eqref{eq:basecamp} is:
	\begin{equation}\label{eq:ffirst}
	\eta \E \left[ (Y - F'(\theta)) \exp \Bigl( \eta \bigl( (\theta - \theta^*) Y + F(\theta^*) - F(\theta) \bigr) \Bigr) \right] .
	\end{equation}
	The second derivative is:
	\begin{equation}\label{eq:fsecond}
	\E \left[ \left( -\eta F''(\theta) + \eta^2 (Y - F'(\theta))^2 \right) 
  \cdot	\exp \Bigl( \eta \bigl( (\theta - \theta^*) Y + F(\theta^*) - F(\theta) \bigr) \Bigr) \right] .
	\end{equation}
	We will also use the standard result \citep{Grunwald07,BarndorffNielsen78} that, since we assume $\theta^* \in \Theta$, 
	\begin{equation}\label{eq:expfacts}
	\begin{aligned}
	\E[Y] = \E_{Y  \sim P_{\theta^*}}[Y] = \mu(\theta^*);  \qquad
	 \text{for all $\theta \in 
		\bar{\Theta}$:} 
	\ F'(\theta) = \mu(\theta); \qquad
	F''(\theta) = \E_{Y  \sim P_{\theta}}(Y- E(Y))^2,
	\end{aligned}
	\end{equation}
	the latter two following because   $F$ is the cumulant generating function.
	
	{\em Part (ii)}.  We use an exact second-order Taylor expansion via the
	Lagrange form of the remainder. We already showed there exist
	$\eta' > 0$ such that, for all $0 < \eta \leq \eta'$,
	all $\theta \in {\Theta}$, $\E[\exp(2\eta (\theta- \theta^*) Y)] <
	\infty$. Fix any such $\eta$.  For some $\theta' \in \left\{ (1 -
	\alpha) \theta + \alpha \theta^* \colon \alpha \in [0, 1] \right\}$,
	the (exact) expansion is:
	\begin{equation} \nonumber
	\begin{multlined}
	\E \left[ e^{-\eta (\xslosstheta)} \right]  
	= 1 + \eta (\theta - \theta^*) \E \left[ Y - F'(\theta^*) \right] \label{eqn:d-cgf-is-mean} 
	- \frac{\eta}{2} (\theta - \theta^*)^2 F''(\theta') 
\cdot	\E \left[ \exp \Bigl( \eta \bigl( (\theta' - \theta^*) Y 
	+ F(\theta^*) - F(\theta')
	\bigr)
	\Bigr) \right]
	 \nonumber \\
	+ \frac{\eta^2}{2} (\theta - \theta^*)^2 	\E \left[ (Y - F'(\theta'))^2 \right.
	\cdot \left. \exp \Bigl( \eta \bigl( (\theta' - \theta^*) Y 
	+ F(\theta^*) - F(\theta')
	\bigr) \Bigr)
	\right] . \nonumber 
	\end{multlined}
	\end{equation}
	Defining $\Delta = \theta' - \theta$, and since $F'(\theta^*) = \E [ Y
	]$ (see \eqref{eq:expfacts}), we see that the central condition is
	equivalent to the inequality:
	\begin{align*}
	\eta \E \left[ (Y - F'(\theta'))^2 e^{\eta \Delta Y} \right] 
	\leq F''(\theta') \E \left[ e^{\eta \Delta Y} \right] .
	\end{align*}
	From Cauchy-Schwarz, to show that the $\eta$-central condition holds it is sufficient to show that 
	\begin{align*}
	\eta \left \| (Y - F'(\theta'))^2 \right\|_{L_2(P)} \left \| e^{\eta \Delta Y} \right\|_{L_2(P)} 
	\leq F''(\theta') \E \left[ e^{\eta \Delta Y} \right] ,
	\end{align*}
	which is equivalent to
	\begin{equation}\label{eqn:eta-expfam-upper-bound}
	\begin{aligned}
	\eta 
	\leq \frac{F''(\theta') \E \left[ e^{\eta \Delta Y} \right]}
	{
		\sqrt{\E \left [ (Y - F'(\theta'))^4 \right] 
			\E \left[e^{2 \eta \Delta Y} \right]}
	} . 
	\end{aligned}
	\end{equation}
	We proceed to lower bound the RHS by lower bounding each of the terms
	in the numerator and upper bounding each of the terms in the
	denominator.  We begin with the numerator. $F'(\theta)$ is bounded by
	\eqref{eq:ineccsi}. Next, by Jensen's inequality,
	$$\E \left[ \exp({\eta \Delta Y}) \right] \geq \exp( \E[\eta 
	\Delta \cdot Y ])\geq \exp(- \eta^{\circ} |\overline\theta -
	\underline\theta| | \mu( \theta^*)| )$$ is lower bounded by a positive
	constant.
	It remains to upper bound the denominator.  Note that the second factor is upper bounded by the constant $\bar{C}$ in \eqref{eq:finite}. 
	The first factor is bounded by a fixed multiple of  
	$\E|Y|^4 + \E[F'(\theta)^4]$. The second term is bounded by \eqref{eq:ineccsi}, so it remains to bound the first term. By assumption $\E[\exp(\eta^{\circ} |Y|)] \leq \bar{C}$ and this implies that $\E|Y^4| \leq a^4 + \bar{C}$ for any $a \geq e$ such that $a^4 \leq  \exp(\eta^{\circ} a)$; such an $a$ clearly exists and only depends on $\eta^{\circ}$. 
	
	We have thus shown that the RHS of \eqref{eqn:eta-expfam-upper-bound} is upper bounded by a quantity that only depends on $\bar{C}, \eta^{\circ}$ and the values of the extrema in \eqref{eq:ineccsi}, which is what we had to show. 
	
	{\em Proof of (iii)}.  We now use the asymptotic form of Taylor's
	theorem. Fix any $\eta > 0$, and pick any $\theta$ close enough to
	$\theta^*$ so that \eqref{eq:basecamp} is finite for all $\theta'$
	in between $\theta$ and $\theta^*$; such a $\theta \neq \theta^*$ must
	exist since for any $\delta > 0$, if
	$|\theta - \theta^*| \leq \delta$, then by assumption
	\eqref{eq:basecamp} must be finite for all
	$\eta \leq \eta^\circ/\delta$.  Evaluating the first and second
	derivative \eqref{eq:ffirst} and \eqref{eq:fsecond} at
	$\theta = \theta^*$ gives:
	\begin{equation*}
	\begin{multlined}
	   \E \left[ e^{-\eta (\xslosstheta)} \right]	= 	1 + \eta (\theta - \theta^*) \E \left[ Y - F'(\theta^*) \right] 
	- \left( \frac{\eta}{2}
	(\theta - \theta^*)^2 F''(\theta^*) 
	- \frac{\eta^2}{2} (\theta - \theta^*)^2 
	 \cdot  \E \left[ (Y -
	F'(\theta^*))^2 \right] \right) \\
	 + h(\theta)(\theta - \theta^*)^2   
	= 1 - \frac{\eta}{2}
	(\theta - \theta^*)^2 F''(\theta^*) 
	+ \frac{\eta^2}{2} (\theta - \theta^*)^2 \E \left[ (Y -
	F'(\theta^*))^2 \right]  + h(\theta) (\theta - \theta^*)^2, 
	\end{multlined}
\end{equation*}
	where $h(\theta)$ is a function satisfying $\lim_{\theta \rightarrow \theta^*} h(\theta) = 0$, where we again used \eqref{eq:expfacts}, i.e.~that $F'(\theta^*) = \E\left[Y\right]$.
	Using further  that $\sigma^2 =  \E \left[ (Y - F'(\theta^*))^2
	\right]$ and $F''(\theta^*) = (\sigma^*)^2$, we find that 
	$\E \left[ e^{-\eta (\xslosstheta)} \right] \leq 1$ iff
	$$
	- \frac{\eta}{2}  (\theta - \theta^*)^2 (\sigma^*)^2
	+ \frac{\eta^2}{2} (\theta - \theta^*)^2 \sigma^2  + h(\theta) (\theta - \theta^*)^2 \leq 0.
	$$
	It follows that for all $\delta > 0$, there is an $\epsilon > 0$ such that for all $\theta \in [\theta^*- \epsilon,\theta^* + \epsilon]$, all $\eta > 0$,
	\begin{align}
	\frac{\eta^2}{2} \sigma^2 \leq \frac{\eta}{2} (\sigma^*)^2  - \delta  
	& \Rightarrow \E \left[ e^{-\eta (\xslosstheta)} \right] \leq 1 \label{eq:goodcondition} \\ \label{eq:badcondition}
	\frac{\eta^2}{2} \sigma^2 \geq \frac{\eta}{2} (\sigma^*)^2  +  \delta  
	& \Rightarrow \E \left[ e^{-\eta (\xslosstheta)} \right] \geq  1\end{align}
	The condition in \eqref{eq:goodcondition} is implied if: 
	$$0 < \eta \leq \frac{(\sigma^*)^2}{\sigma^2} -   \frac{2 \delta}{\eta \sigma^2}.$$
	Setting $C = 4 \sigma^2/(\sigma^*)^4$ and
	$\eta_{\delta} = (1- C {\delta}) (\sigma^*)^2/\sigma^2$ we find that
	for any $\delta < (\sigma^*)^4/(8 \sigma^2)$, we have
	$1- C {\delta} \geq 1/2$ and thus $\eta_{\delta} > 0$ so that in
	particular the premise in \eqref{eq:goodcondition} is satisfied for
	$\eta_{\delta}$.  Thus, for all small enough $\delta$, both the
	premise and the conclusion in \eqref{eq:goodcondition} hold for
	$\eta_{\delta} > 0$; since
	$\lim_{\delta \downarrow 0} \eta_{\delta} = (\sigma^*)^2/\sigma^2$, it
	follows that there is an increasing sequence
	$\eta_{(1)}, \eta_{(2)}, \ldots$ converging to $(\sigma^*)^2/\sigma^2$
	such that for each $\eta_{(j)}$, there is $\epsilon_{(j)} > 0$ such
	that for all
	$\theta \in [\theta^*- \epsilon_{(j)}, \theta^*+\epsilon_{(j)}]$,
	$\E \left[ e^{-\eta_{(j)} (\xslosstheta)} \right] \leq 1$. It follows
	that the $\lim \sup$ in \eqref{eq:limitvariance} is at least
	$(\sigma^*)^2/\sigma^2$. A similar argument (details omitted) using
	\eqref{eq:badcondition} shows that the $\lim \sup$ is at most this
	value; the result follows.
\subsection{Proof  of Proposition~\ref{prop:entroboundb}}
For arbitrary conditional densities $p'(y \mid x)$ with corresponding distribution $P'\mid X$ for which
	\begin{equation}\label{eq:daling}
	\E_{P'}[Y |X] = g^{-1}(\langle \beta,X),
	\end{equation} 
	and 
	densities $p_{\fopt} = p_{\beta^*}$ and $p_\beta$ with $\beta^*,\beta \in \cB$, we can write:
	\begin{equation*}
	\begin{aligned}
	\E_{X\sim P} \E_{Y \sim P'\mid X}\left[ \log  \frac{p_{\beta^*}(Y \mid X)}{p_{\beta}(Y | X)}\right] 
	&=   \E \E\left[ (\theta_X(\beta^*) - \theta_X(\beta) )  Y -  \log 
	\frac{F(\theta_X(\beta^*))}{F(\theta_X(\beta))}   \mid X \right]\\
	 &= \E_{X \sim P} 
	\left[ (\theta_X(\beta^*) - \theta_X(\beta) )g^{-1}(\langle \beta, X\rfloor_d \rangle 
 -  \log F(\theta_X(\beta^*)) + \log F(\theta_X(\beta))   \mid X \right],
	\end{aligned}
	\end{equation*}
	where the latter equation follows by \eqref{eq:daling}. The result now
	follows because \eqref{eq:daling} both holds for the `true' $P$ and
	for $P_{\fopt}$.
\subsection{Proof of Proposition~\ref{prop:tenerife}}
	The fact that under the three imposed conditions the
	$\bar\eta$-central condition holds for some $\bar\eta > 0$ is a simple consequence of
	Theorem~\ref{thm:expfam-exptails}: Condition 1 implies that there
	is some compact $\Theta$ such that for all $x \in \cX$, $\beta \in
	\cB$, $\theta_x(\beta) \in \Theta$. Condition 3 then ensures that
	$\theta_x(\beta)$ lies in the interior of this $\Theta$. And Condition
	2 implies that $\bar\eta$ in Theorem~\ref{thm:expfam-exptails}
	can be chosen uniformly for all $x \in \cX$.
\section{EXCESS RISK AND KL DIVERGENCE INSTEAD OF GENERALIZED HELLINGER DISTANCE}
\label{sec:witness}
The misspecification metric/generalized Hellinger distance 
$\dhel_{\bar\eta}$ appearing in Theorem~\ref{thm:metric} is rather weak (it is `easy' for two distributions
to be close) and lacks a clear interpretation for general,
non-logarithmic loss functions. Motivated by these facts, GM study in
depth under what additional conditions the (square of this) metric can
be replaced by a stronger and more readily interpretable divergence
measure. They come up with a
new, surprisingly weak condition, the {\em witness condition}, under which $\dhel_{\bar\eta}$ can be
replaced by the {\em excess risk\/} ${\bf E}_P[L_f]$, which is the
additional risk incurred by $f$ as compared to the optimal $f^*$. For
example, with the squared error loss, this is the additional mean
square error of $f$ compared to $f^*$; and with (conditional)
log-loss, it is the well-known {\em generalized KL divergence\/} ${\bf
  E}_{X,Y \sim P}[\log \frac{p_{f^*}(Y \mid X)}{p_f(Y|X)}]$, coinciding with
  standard KL divergence if the model is correctly specified. Bounding
  the excess risk is a standard goal in statistical learning theory;
  see for example \citep{bartlett2005local,erven2015fast}.
  
The following definition appears (with substantial explanation
including the reason for its name) as Definition 12 in GM:
\begin{definition}[Empirical Witness of Badness] \label{def:witness}
	We say that $\smtuple$ satisfies the $(u,c)$-\emph{empirical witness of badness condition} (or \emph{witness condition}) 
	for constants $u > 0$ and $c \in (0,1]$ if for all $f \in \cF$
	\begin{align*}
	\E \left[ (\xslosslong{f}) \cdot \ind{\xslosslong{f} \leq u} \right] \geq c \xsrisklong{f} .
	\end{align*}
	More generally, for a function $\tau: \reals^+ \to [1,\infty)$
	and constant $c \in (0,1)$ 
	we say that $\smtuple$ satisfies the
	\emph{$(\tau,c)$-witness  condition} if 
	for all $f \in \cF$, $\xsrisklong{f} < \infty$ and
	\begin{align*} 
	\E \left[ (\xslosslong{f}) \cdot \ind{\xslosslong{f} \leq \tau(\xsrisklong{f})} \right] \geq c \xsrisklong{f} .
	\end{align*}
\end{definition}
It turns out that the $(\tau,c)$-witness condition holds in many
practical situations, including our GLM-under-misspecification
setting. Before elaborating on this, let us review (a special case of)
Theorem~12 of GM, which is the analogue of Theorem~\ref{thm:metric}
but with the misspecification metric replaced by the excess risk.

First, let, for arbitrary $0 < \eta < \bar\eta$,  
	$c_u \coloneqq \frac{1}{c} \frac{\eta u + 1}{1 - \frac{\eta}{\bar{\eta}}}$.  Note that for large $u$, $c_u$ is approximately linear in $u/c$.
\begin{theorem}{\bf [Specialization of Theorem 12 of GM]}
	Consider a learning problem $\smtuple$.  Suppose that the $\bar\eta$-strong central condition
	holds.  If the $(u,c)$-witness condition holds, then for any
	$\eta \in (0,\bar\eta)$,
	\begin{align}\label{eq:same}
	{\bf E}_{Z^n \sim P} \E_{\rv{f} \sim \dol_n}\left[\xsrisk{f}\right] 
	\leq   
	c_u \cdot {\bf E}_{Z^n \sim P} \left[
          \rsc_{n,\eta} \left( \Pi_0 \right) \right],
	\end{align}
	with $c_u$ as above. If instead the $(\tau,c)$-witness
        condition holds for some nonincreasing function $\tau$ as
        above, then for any $\lambda > 0$,
	\begin{equation}
	{\bf E}_{Z^n \sim P} \E_{\rv{f} \sim \dol_n} \left[\xsrisk{f} \right] 
	\ \leq  \ 
	\lambda + c_{\tau(\lambda)} \cdot {\bf E}_{Z^n \sim P} \left[ \rsc_{n,\eta} \left( \Pi_0\right)\right]. \nonumber
	\end{equation}
\end{theorem}
The actual theorem given by GM generalizes this to an in-probability
statement for general (not just generalized Bayesian) learning
methods.  If the $(u,c)$-witness condition holds, then, as is obvious
from (\ref{eq:same}) and Theorem~\ref{thm:metric}, the same rates can
be obtained for the excess risk as for the squared misspecification
metric.  For the $(\tau,c)$-witness condition things are a bit more
complicated; the following lemma (Lemma 16 of GM) says that, under an
exponential tail condition, $(\tau,c)$-witness holds for a
sufficiently `nice' function $\tau$, for which we loose at most a
logarithmic factor:
\begin{lemma} \label{lem:kl-hell-exp-tails} Define
	$M_{\kappa} \coloneqq \sup_{f \in \cF} \E \left[ e^{\kappa \xsloss{f}} \right]$ and assume
	that the excess loss $L_f$ has a uniformly exponential upper tail, i.e.\ $M_{\kappa} < \infty$.  Then, for the map
	$\tau: x \mapsto 1 \opmax \kappa^{-1} {\log \frac{2 M_\kappa }{\kappa x}} = O(1 \opmax \log (1/x))$, the $(\tau,c)$-witness condition holds with $c = \nicefrac{1}{2}$. 
\end{lemma}
As an immediate consequence of this lemma, 
GM's theorem above gives that for any $\eta \in (0,\bar\eta)$, (using $\lambda=
1/n$), there is $C_{\eta} < \infty$ such that
\begin{align} \label{eqn:kl-hell-exp-tails-2} 
{\bf E}_{Z^n \sim P} \E_{\rv{f} \sim \dol_n} \left[ \xsrisk{\rv{f}} \right] \leq 
\frac{1}{n}
+ C_{\eta} \cdot (\log n) \cdot {\bf E}_{Z^n \sim P} \left[\rsc_{\eta,n} \left( \fopt \pipes  \dolest \right) \right],
\end{align}
so our excess risk bound is only a log factor worse than the bound
that can be obtained for the squared misspecification metric in
Theorem~\ref{thm:metric}. We now apply this to the misspecified GLM setting:
\paragraph{Generalized Linear Models and Witness}
	Recall that the central condition holds for generalized linear
        models under the three assumptions made in
        Proposition~\ref{prop:tenerife}.  Let $\loss_\beta \coloneqq
        \loss_\beta(X,Y) = -\log p_\beta(Y \mid X)$ be the loss of
        action $\beta \in \cB$ on random outcome $(X,Y) \sim P$, and
        let $\beta^*$ denote the risk minimizer over $\cB$.  The first
        two assumptions taken together imply, via \eqref{eq:ineccsi},
        that there is a $\kappa > 0$ such that
	\begin{align*}
	\sup_{\beta \in B} \E_{X,Y \sim P}\left[ 
	e^{\kappa (\loss_\beta - \loss_{\beta^*})} \right] 
	&\leq 
	\sup_{\beta \in \cB, x \in \cX} \E_{Y \sim P \mid X=x}\left[ 
	e^{\kappa (\loss_\beta - \loss_{\beta^*})}
	\right] \\ 
	&= \sup_{\beta \in \cB, x \in \cX} 
	\left( \frac{F_{\theta_x(\beta)}}{F_{\theta_x(\beta^*)}}\right)^{\kappa} \cdot \E_{Y \sim P \mid X=x}\left[ 
	e^{\kappa | Y| }\right]
	< \infty.
	\end{align*}
	The conditions of Lemma~\ref{lem:kl-hell-exp-tails} are thus
        satisfied, and so the $(\tau,c)$-witness condition holds for
        the $\tau$ and $c$ in that lemma. From
        \eqref{eqn:kl-hell-exp-tails-2} we now see that we get an
        $O((\log n)^2/n)$ bound on the expected excess risk, which is
        equal to the parametric (minimax) rate up to a $(\log n)^2$
        factor. Thus, fast learning rates in terms of excess risks and
        KL divergence under misspecification with GLMs are possible
        under the conditions of Proposition~\ref{prop:tenerife}.

\section{LEARNING RATE $> 1$ FOR MISSPECIFIED MODELS}\label{app:learningrate>1}

In what follows we give an example of a misspecified setting, where the best performance is achieved with the learning rate $\eta>1$. Consider a model $\{P_\beta,\beta\in[0.2,0.8]\},$ where $P_\beta$ is a Bernoulli distribution with $\PP_\beta(Y=1)=\beta$. Let the data $Y_1,\dots,Y_n$ be sampled i.i.d.\ from $P_0$, i.e.\ $Y_i=0$ for all $i=1,\dots,n$. In this case the log-likelihood function is given by
\[
\log p(Y_1,\dots,Y_n\given \beta)=n\log(1-\beta).
\]
Observe that in this setting $\beta^\star=0.2$. Now assume that the model is correct and data $Y_1',\dots,Y_n'$ is sampled i.i.d.\ from $P_{
\beta}$ with $\beta=0.2$. Then the log-likelihood is
\[
\log p(Y_1',\dots,Y_n'\given \beta=0.2)\approx 0.2n\log0.2+0.8n\log0.8\ll n\log0.8=\log p(Y_1,\dots,Y_n\given \beta=0.2).
\]
Thus, the data are more informative about the best distribution than they would be if the model were correct. Therefore, we can afford to learn ‘faster’: let the data be more important and the (regularizing) prior be less important. This is realized by taking $\eta >> 1$.


\section{MCMC SAMPLING
}\label{app:impl}

\subsection{The $\eta$-generalized Bayesian lasso}
Here, following \citet{Park} we consider a slightly more general version of the regression problem:
\[
Y=\mu+  X  \beta+\eps,
\]
where $\mu\in\RR^n$ is the overall mean, $\beta\in\RR^p$ is the vector of parameters of interest, $y\in\RR^n$, $X\in\RR^{n\times p}$, and $\eps\sim N(0,\sigma^2 I_n)$ is a noise vector. For a given shrinkage parameter $\lambda>0$ the Bayesian lasso of \citet{Park} can be represented as follows.
\begin{align}\label{eq:hiermodelbayesianlassopark}
 {Y} \rvert \mu,  {X},  {\beta}, \sigma^2 &\sim N(\mu +  {X\beta}, \sigma^2 {I_n}) \, , \\ \nonumber
 {\beta} \rvert \tau_1^2, \ldots, \tau_p^2, \sigma^2 &\sim N( {0}, \sigma^2  {D}_\tau), \, \, \, \,  {D}_\tau = \diag(\tau_1^2, \ldots, \tau_p^2) \, , \\ \nonumber
\tau_1^2, \ldots, \tau_p^2 &\sim \prod_{j=1}^p \frac{\lambda^2}{2} e^{-\lambda^2\tau_j^2 /2} d\tau_j^2, \, \, \, \, \tau_1^2, \ldots, \tau_p^2 > 0 \, , \\ \nonumber
\sigma^2 &\sim \pi(\sigma^2) \, d\sigma^2 \, .
\end{align}

In this model formulation the $\mu$ on which the outcome variables $Y$ depend, is the overall mean, from which ${X\beta}$ are deviations. The parameter $\mu$ can be given a flat prior and subsequently integrated out, as we do in the coming sections.

We will use the typical inverse gamma prior distribution on $\sigma^2$, i.e.\ for $\sigma^2>0$
\[
\pi(\sigma^2)=\frac{\gamma^\alpha}{\Gamma(\alpha)}\sigma^{-2\alpha-2}e^{-\gamma/\sigma^2},
\]
where $\alpha,\gamma>0$ are hyperparameters. With the hierarchy of \eqref{eq:hiermodelbayesianlassopark}  the joint density for the posterior with the likelihood to the power $\eta$ becomes

\begin{multline}\label{eq:joint generalized density}
(f({Y} \rvert \mu,  {\beta}, \sigma^2))^\eta \, \pi(\sigma^2) \, \pi(\mu) \prod_{j=1}^p \, \pi( {\beta}_j \rvert \tau_j^2, \sigma^2 ) \, \pi(\tau_j^2) = \\ 
=\left( \frac{1}{(2\pi\sigma^2)^{n/2}} \,  e^{\frac{1}{2\sigma^2} ( {Y} -\mu  {1}_n -  {X \beta})^T  ( {Y} -\mu  {1}_n -  {X \beta})  }   \right) ^\eta 
\frac{\gamma^\alpha}{\Gamma(\alpha)}\sigma^{-2\alpha-2}e^{-\frac{\gamma}{\sigma^2}} \prod_{j=1}^p \frac{1}{(2\sigma^2\tau_j^2)^{1/2}} \, e^{- \frac{1}{2\sigma^2\tau_j^2}  {\beta}_j^2} \frac{\lambda^2}{2} \, e^{- \lambda^2\tau^2_j /2} \,.
\end{multline}

Let $ {\tilde{Y}}$ be $ {Y - \overline{Y}}$. If we integrate out $\mu$, the joint density marginal over $\mu$ is proportional to

\begin{align}\label{eq:joint generalized density marginalized over mu}
\sigma^{-\eta(n-1)} \,  e^{-\frac{\eta}{2\sigma^2} ( {\tilde{Y}-X\beta})^T  ( {\tilde{Y}-X\beta})  } \,
\sigma^{-2\alpha-2} \, e^{-\frac{\gamma}{\sigma^2}} \prod_{j=1}^p \frac{1}{(\sigma^2\tau_j^2)^{1/2}} \, e^{- \frac{1}{2\sigma^2\tau_j^2}  {\beta}_j^2} \, e^{- \lambda^2\tau^2_j /2} .
\end{align}

First, observe that the full conditional for $  \beta$ is multivariate normal: the exponent terms involving $ {\beta}$ in \eqref{eq:joint generalized density marginalized over mu} are 

\[
-\frac{\eta}{2\sigma^2} ( {\tilde{Y}-X\beta})^T  ( {\tilde{Y}-X\beta}) - \frac{1}{2\sigma^2}  {\beta}^T  {D_\tau}^{-1}  {\beta}= -\frac{1}{2\sigma^2} \left\lbrace ( {\beta}^T (\eta  {X}^T {X} +  {D_\tau}^{-1} ) {\beta} - 2\eta {\tilde{Y}X\beta} + \eta {\tilde{Y}}^T  {\tilde{Y}} ) \right\rbrace.
\]

If we now write $M_\tau = (\eta  {X}^T {X} +  {D_\tau}^{-1})^{-1}$ and complete the square, we arrive at
\begin{align*}
&-\frac{1}{2\sigma^2} \left\lbrace ( {\beta} - \eta M_\tau  {X}^T  {\tilde{Y}})^T {M_\tau^{-1}} \, ( {\beta} - \eta {M_\tau} {X}^T  {\tilde{Y}}) +  {\tilde{Y}}^T (\eta  {I}_n - \eta^2  {X}^{-1}M_\tau {X}^T)  {\tilde{Y}} \right\rbrace .
\end{align*}

Accordingly we can see that $ {\beta}$ is conditionally multivariate normal with mean $\eta {M_\tau} {X}^T  {\tilde{Y}}$ and variance $\sigma^2 {M_\tau}$.\\

The terms in \eqref{eq:joint generalized density marginalized over mu} that involve $\sigma^2$ are:

\[
(\sigma^2)^{ \{ -\eta(n-1)/2 - p/2 - \alpha - 1 \} }
\exp \Big\{ -\frac{\eta}{2\sigma^2}( {\tilde{Y}-X\beta})^T  ( {\tilde{Y}-X\beta}) - \frac{1}{2\sigma^2} {\beta}^T  {D_\tau}^{-1}  {\beta} - \frac{\gamma}{\sigma^2}  \Big\} .
\]

We can conclude that $\sigma^2$ is conditionally inverse gamma with shape parameter \linebreak $\displaystyle \eta \, \frac{n-1}{2} + \frac{p}{2} + \alpha$ and scale parameter $\displaystyle \frac{\eta}{2}( {\tilde{Y}-X\beta})^T  ( {\tilde{Y}-X\beta}) +  {\beta}^T  {D_\tau}^{-1}  {\beta} /2 + \gamma$. \\
\\
Since $\tau_j^2$ is not involved in the likelihood, we need not modify the implementation of it and follow \citet{Park}:
\[
\frac{1}{\tau^2_j} \sim \text{IG}\left(\displaystyle \sqrt{{\lambda^2\sigma^2}/{\beta_j^2}}, \, \lambda^2\right). 
\]
Summarizing, we can implement a Gibbs sampler with the following distributions:

\begin{align} \label{Blasso:condbetagen}
 {\beta} &\sim \text{N} \left( \eta(\eta {X}^T {X} +  {D_\tau}^{-1})^{-1} {X}^T {\tilde{Y}}, \, \sigma^2(\eta {X}^T {X} +  {D_\tau}^{-1})^{-1} \right) \, , \\ \label{Blasso:condsigmagen}
\sigma^2 & \sim \text{Inv-Gamma} \big( \frac{\eta}{2}(n-1) + p/2 + \alpha,\, \frac{\eta}{2}( {\tilde{Y}-X\beta})^T  ( {\tilde{Y}-X\beta}) +  {\beta}^T  {D_\tau}^{-1}  {\beta}/2 + \gamma \big) \, , \\ \label{Blasso:condtausqgen}
\frac{1}{\tau^2_j} &\sim \text{IG}\left(\displaystyle \sqrt{{\lambda^2\sigma^2}/{\beta_j^2}}, \, \lambda^2\right) \, .
\end{align}
There are several ways to deal with the shrinkage parameter $\lambda$. We follow the hierarchical Bayesian approach and place a hyperprior on the parameter.  In our implementation we provide three ways to do so: a point mass (resulting in a fixed $\lambda$), a gamma prior on $\lambda^2$ following \citet{Park} and a beta prior following \citet{deloscampos}, details about the implementation of the latter two priors can be found in those papers respectively.

\subsection{The $\eta$-generalized Bayesian logistic regression}\label{ap:logreg}
We follow the construction of the P{\'o}lya--Gamma latent variable scheme for constructing a Bayesian estimator in the logistic regression context described in \cite{polson2013bayesian}. 

First, for $b>0$ consider the density function of a P{\'o}lya-Gamma random variable $PG(b,0)$
\[
p(x\given b,0)=\frac{2^{b-1}}{\Gamma(b)}\sum_{n=1}^\infty (-1)^n\frac{\Gamma(n+b)}{\Gamma(n+1)}\frac{(2n+b)}{\sqrt{2\pi x^3}}e^{-\frac{(2n+b)^2}{8x}}.
\] 
The general class $PG(b,c)$ ($b,c>0$) is defined through an exponential tilting of the $PG(b,0)$ and has the density function 
\[
p(x\given b,c)=\frac{e^{-\frac{c^2x}2}p(x|b,0)}{\EE e^{-\frac{c^2\omega}2}},
\]
where $\omega\sim PG(b,0)$.

To derive our Gibbs sampler we use the following result from \cite{polson2013bayesian}.

\begin{thm}
	\label{PG_thm}
	Let $p_{b,0}(\omega)$ denote the density of $PG(b,0)$. Then for all $a\in \RR$
	\[
	\frac{(e^\psi)^a}{(1+e^\psi)^b}=2^{-b}e^{\kappa \psi}\int_0^\infty e^{-\omega\psi^2/2}p_{b,0}(\omega)d\omega,
	\]
	where $\kappa=a-b/2.$
\end{thm}

According to Theorem \ref{PG_thm} the likelihood contribution of the observation $i$ taken to the power $\eta$ can be written as
\[
L_{i,\eta}(\beta)=\left[\frac{(e^{X_i^T\beta)^{y_i }}}{1+e^{X_i^T\beta}}\right]^\eta\propto e^{\eta\kappa_iX_i^T\beta}\int_0^\infty e^{-\omega_i\frac{(X_i^T\beta)^2}2}p(\omega_i\given \eta,0),
\]
where $\kappa_i\coloneqq y_i-1/2$ and $p(\omega_i\given \eta,0)$ is the density function of $PG(\eta,0)$. 

Let 
\begin{gather*}
X\coloneqq (X_1,\dots, X_n)^T, \quad Y\coloneqq (Y_1,\dots,Y_n)^T, \quad \kappa\coloneqq (\kappa_1,\dots,\kappa_n)^T, \\
\omega\coloneqq (\omega_1,\dots,\omega_n)^T,  \quad \Omega\coloneqq \diag(\omega_1,\dots,\omega_n).
\end{gather*}
Also, denote the density of the prior on $\beta$ by $\pi(\beta)$. Then the conditional posterior of $\beta$ given $\omega$ is
\[
p(\beta\given \omega,Y)\propto \pi(\beta) \prod_{i=1}^n L_{i,\eta}(\beta\given\omega_i)=\pi(\beta)\prod_{i=1}^ne^{\eta\kappa_iX_i^T\beta-\omega_i\frac{(X_i^T\beta)^2}2}\propto \pi(\beta)e^{-\frac12 (z-X\beta)^T\Omega(z-X\beta)},
\]
where $z\coloneqq \eta(\frac{\kappa_1}{\omega_1},\dots,\frac{\kappa_n}{\omega_n})$. Observe that the likelihood part is conditionally Gaussian in $\beta$. Since the prior on $\beta$ is Gaussian, a simple linear-model calculation leads to the following Gibbs sampler. To sample from the the $\eta$-generalized posterior one has to iterate these two steps
\begin{align}
\omega_i\given \beta\sim& PG(\eta, X_i^T\beta),\\
\beta\given Y,\omega\sim& \cN(m_\omega, V_\omega), 
\end{align}
where
\begin{align*}
V_\omega\coloneqq &(X^T\Omega X+B^{-1})^{-1}, \\
m_\omega\coloneqq &V_\omega(\eta X^T\kappa+B^{-1}b).
\end{align*}

To sample from the P{\'olya}-Gamma distribution $PG(b,c)$ we adopt a method from \citep{PG2014}, which is based on the following representation result. According to \cite{polson2013bayesian} a random variable $\omega\sim PG(b,c)$ admits the following representation
\[
\omega \eqd \sum_{n=0}^\infty \frac{g_n}{d_n},
\]
where $g_n\sim Ga(b,1)$ are independent Gamma distributed random variables, and 
\[
d_n\coloneqq 2\pi^2(n+\frac12)^2+2c^2.
\]

Therefore, we approximate the PG random variable by a truncated sum of weighted Gamma random variables. \citep{PG2014} shows that the approximation method performs well with the truncation level $N=300$. Furthermore, we performed our own comparison of the sampler with the STAN implementation for Bayesian logistic regression, which showed no difference between the methods (for $\eta=1$).

\subsection{The Safe-Bayesian Algorithms}\label{ap:safebayesimpl}

%

 The version of the Safe-Bayesian algorithm we are using for the experiments is called \emph{R-log-SafeBayes}, more details and other versions can be found in \citet{grunwald2017inconsistency}. The $\hat{\eta}$ is chosen  from a grid of learning rates $\eta$ that minimizes the \emph{cumulative Posterior-Expected Posterior-Randomized log-loss}:
\begin{align*}
\sum_{i=1}^n \, \E_{\beta, \sigma^2 \sim \Pi \rvert z^{i-1}, \eta}  \left[-\log f(Y_i \rvert X_i, \beta, \sigma^2) \right].
\end{align*}
Minimizing this comes down to minimizing
\begin{align*}
\sum_{i=1}^{n-1} \, \, \textsc{av} \left[ \frac{1}{2} \log 2 \pi \sigma^2_{i,\eta} + \frac{1}{2} \frac{(Y_{i+1} - X_{i+1}\beta_{i, \eta})^2}{\sigma^2_{i,\eta}} \right].
\end{align*}
The loss between the brackets is averaged over many draws of $(\beta_{i, \eta}, \sigma^2_{i,\eta})$ from the posterior, where $\beta_{i, \eta}$ (or $\sigma^2_{i,\eta}$) denotes one random draw from the conditional $\eta$-generalized posterior based on data points $z^{i}$. For the sake of completeness we present the algorithm below.\\
\begin{algorithm}
	\SetKwInOut{Input}{Input}\SetKwInOut{Output}{Output}
	\Input{data $z_1, \ldots z_n$, model $\mathcal{M} = \{ f(\cdot \rvert \theta) \rvert \theta \in \Theta \}$, prior $\Pi$ on $\Theta$, step-size $\mathcal{K}_{\textsc{step}}$,  max. exponent $\mathcal{K}_{\textsc{max}}$, loss function $\ell_\theta(z)$ } 
	\Output{Learning rate $\hat{\eta}$} 
	
	$\mathcal{S}_n \coloneqq \{1, 2^{-\mathcal{K}_{\text{STEP}}}, 2^{-2\mathcal{K}_{\text{STEP}}}, 2^{-3\mathcal{K}_{\text{STEP}}}, \ldots, 2^{-\mathcal{K}_{\text{MAX}}}, \} $ \;
	\For{all $\eta \in \mathcal{S}_n$}{
		$s_\eta \coloneqq 0$ \;
		\For{$i = 1 \ldots n$}{
			Determine generalized posterior $\Pi ( \cdot \rvert z^{i-1}, \eta )$ of Bayes with learning rate $\eta$. \\ Calculate posterior-expected posterior-randomized loss of predicting actual next outcome:
			\begin{align}
			r \coloneqq \ell_{\Pi \rvert z^{i-1}, \eta} (z_i) = \E_{\theta \sim \Pi \rvert z^{i-1}, \eta} \left[ \ell_\theta (z_i) \right] 
			\end{align} \\
			$s_\eta \coloneqq s_\eta + r$ \;
		}
	}
	Choose $\hat{\eta} \coloneqq \argmin_{\eta \in \mathcal{S}_n} \{ s_\eta \}$ (if min achieved for several $\eta \in \mathcal{S}_n$, pick largest) \;
	\caption{\emph{The R-Safe-Bayesian algorithm}} \label{alg:r-safebayesian}
\end{algorithm}

\section{DETAILS FOR THE EXPERIMENTS AND FIGURES}\label{ap:detailsfigures}

Below we present the results of additional simulation experiments for Section~\ref{sec:lassoandhs} (Appendix~\ref{ap:simulations}) and the description of experiments with real-world data (Appendix~\ref{ap:realworld}). We also give details for Figure~\ref{fig:predvar} in Appendix~\ref{ap:delete}.

\subsection{Additional Figures for Section~\ref{sec:lassoandhs}}
\label{ap:simulations}

\begin{figure}
	\centerline{\includegraphics[width=0.5\textwidth]{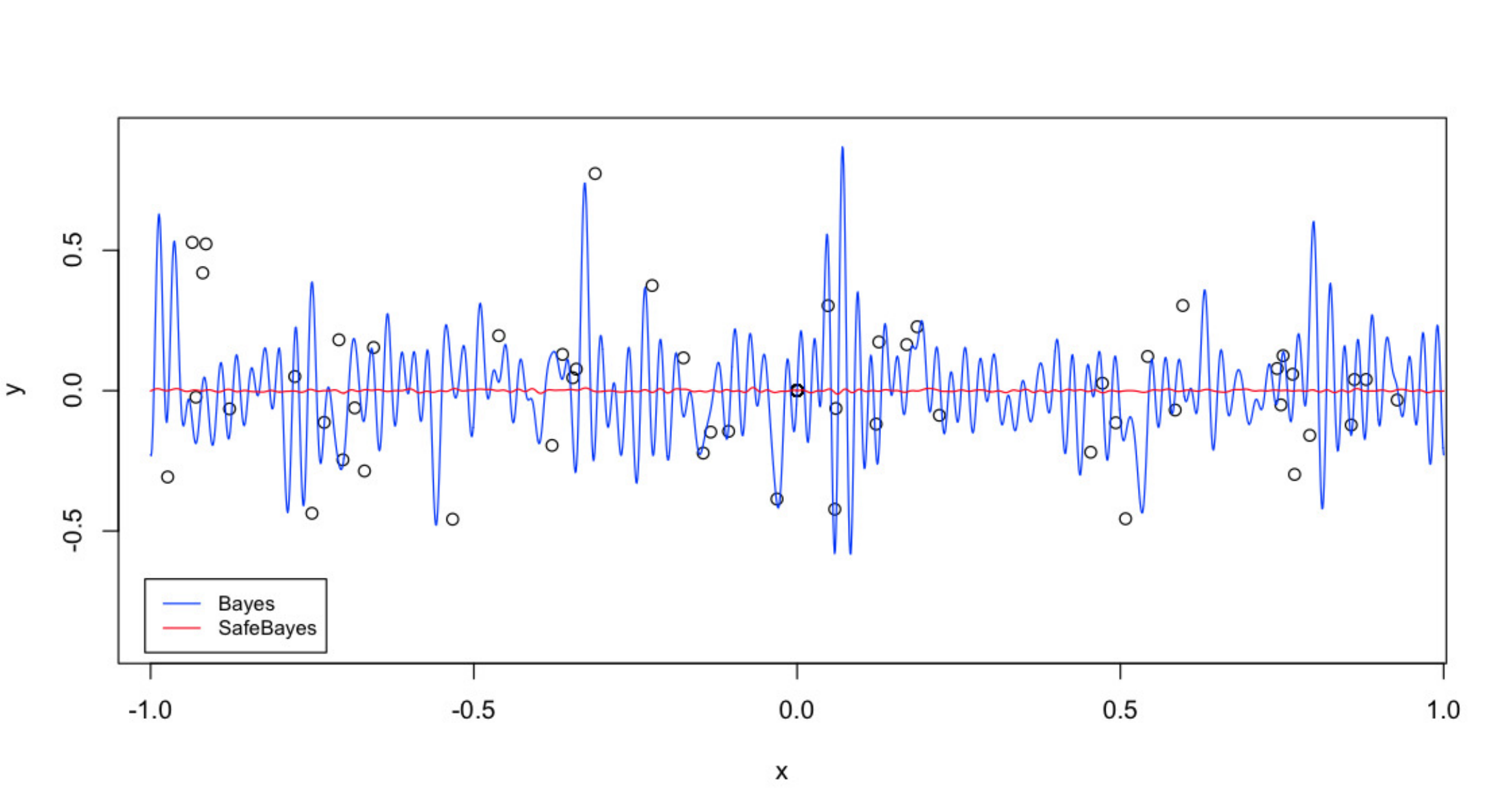}}
	\caption{\emph{Prediction of standard Bayesian lasso (blue) and Safe-Bayesian lasso (red, $\eta = 0.5$) with $n=200$, $p=100$.}\label{fig:b vs sb}}
\end{figure}

\begin{figure}
	\centerline{\includegraphics[width=0.5\textwidth]{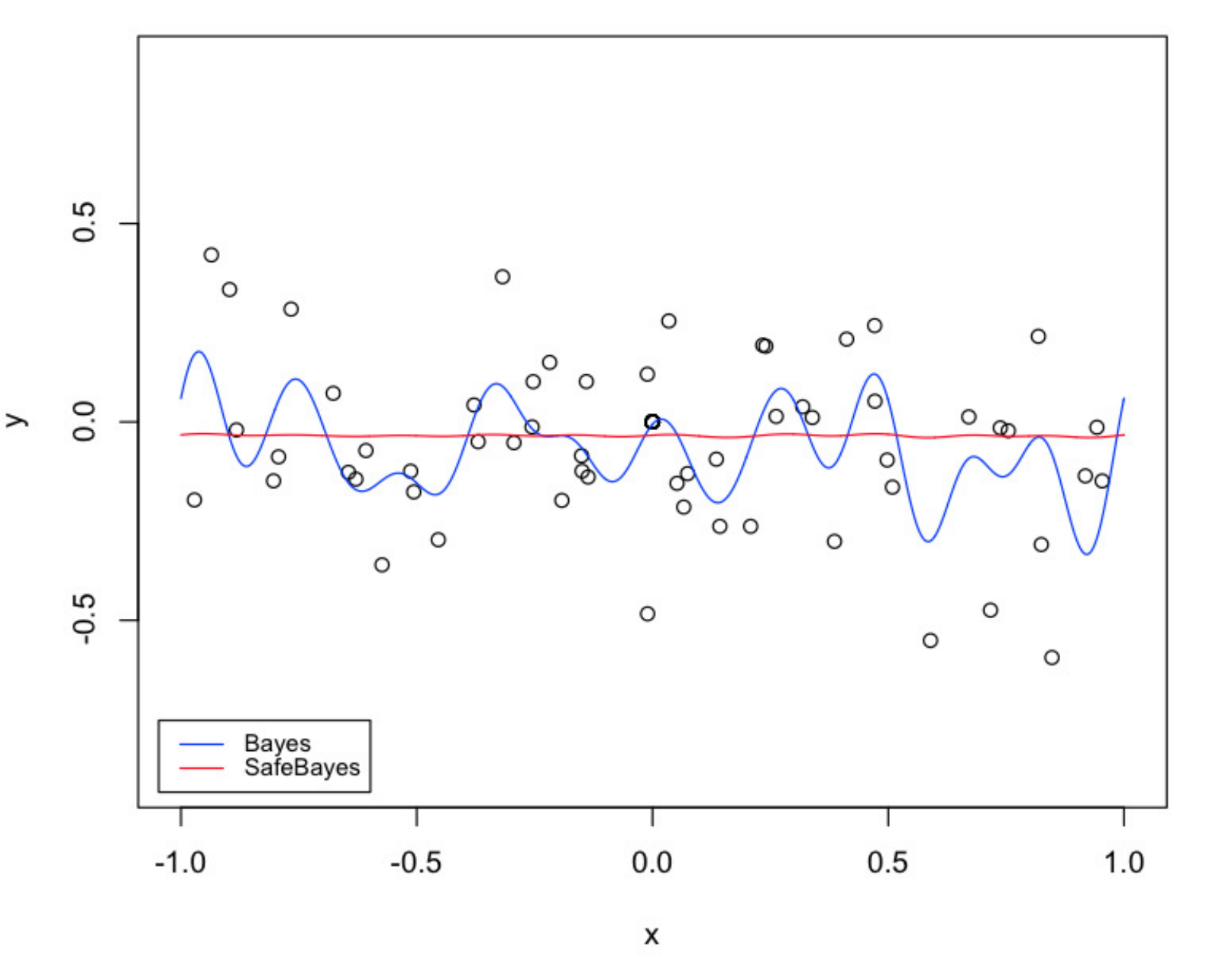}}
	\caption{\emph{Prediction of standard Bayesian lasso (blue) and Safe-Bayesian lasso (red, $\eta = 0.5$) with $n=200$, $p=12$.}\label{fig:1}}
\end{figure}

Consider the regression context described in Section~\ref{sec:lassoandhs}. Here, we explore different choices of the number of Fourier basis functions, showing that regardless of the choice Safe-Baysian lasso outperforms its standard counterpart. In Figures~\ref{fig:b vs sb} and \ref{fig:1} we see conditional expectations $\E \left[ Y \mid X \right]$ according to the posteriors of the standard Bayesian lasso (blue) and the Safe-Bayesian lasso (red, $\hat\eta = 0.5$) for the \emph{wrong-model} experiment described in Section \ref{sec:lassoandhs}, with $100$ data points. We take  $201$ and $25$ Fourier basis functions respectively.  
%
%
%

\begin{figure}
	\centerline{\includegraphics[width=0.45\textwidth]{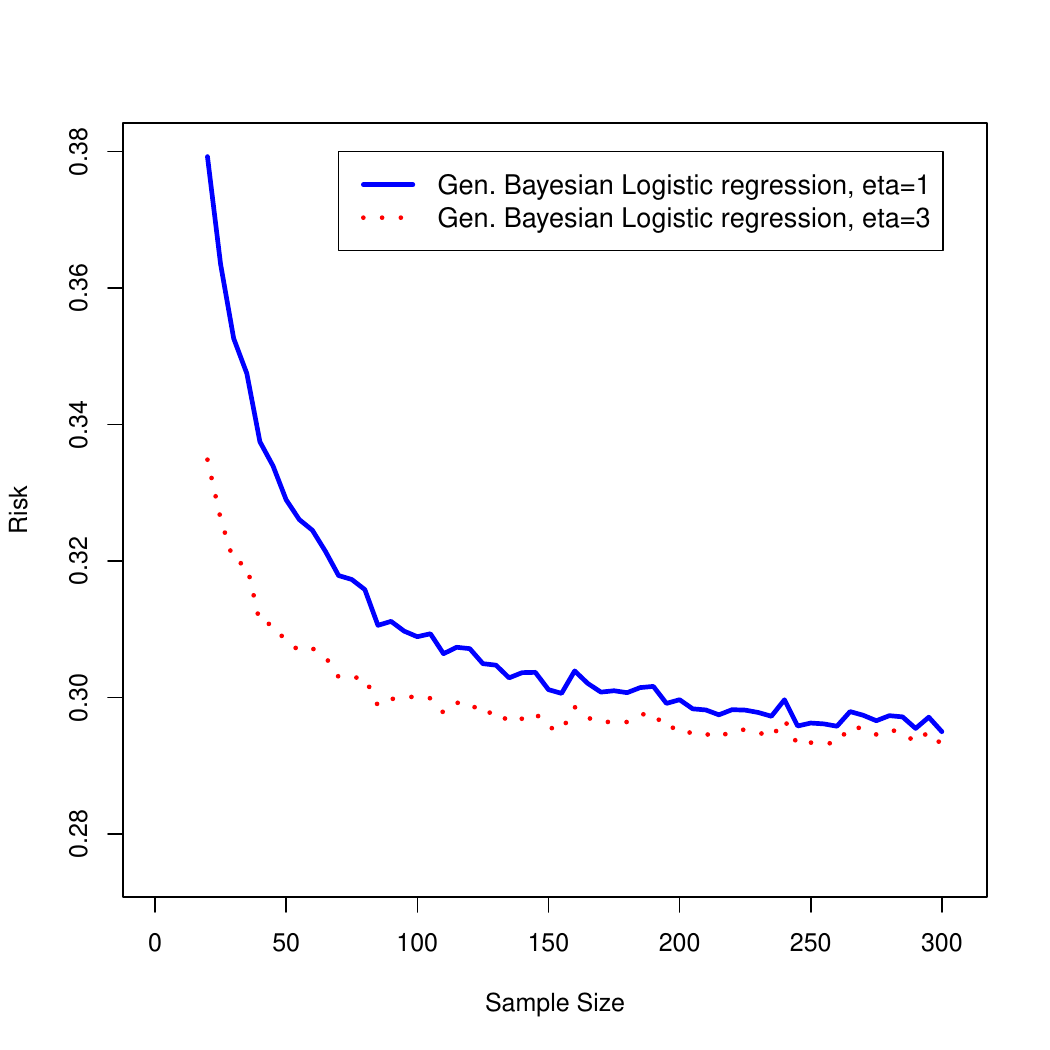}}
	\caption{\emph{Simulated logistic risk as a function of the sample size for the \emph{correct-model} experiments described in Section \ref{sec:logisticregression} according to the posterior predictive distribution of standard Bayesian logistic regression ($\eta=1$), and generalized Bayes ($\eta=3$).} \label{fig:logregriskwell}}
\end{figure}

Now we consider logistic regression setting and show that even for some well-specified problems it is beneficial to choose $\eta\neq 1$. In Figure~\ref{fig:logregriskwell} we see a comparison of the log-risk for $\eta=1$ and $\eta=3$ in the well-specified logistic regression case (described in Section \ref{sec:logisticregression}). Here $p=1$ and $\beta=4$.
 
\subsection{Real-world data}
\label{ap:realworld}
\paragraph{Seattle Weather Data}\label{ap:seattle}
The $\texttt{R}$-package $\texttt{weatherData}$ \citep{weatherData} loads weather data available online from $\texttt{www.wunderground.com}$.  Besides data from many thousands of personal weather stations and government agencies, the website provides access to data from Automated Surface Observation Systems (ASOS) stations located at airports in the US, owned and maintained by the Federal Aviation Administration. Among them is a weather station at Seattle Tacoma International Airport, Washington (WMO ID $72793$). From this station we collected the data for this experiment.

The training data are the maximum temperatures for each day of the year 2011 at Seattle airport. We divided the data randomly in a training set (300 measurements) and a test set (65 measurements). First, we sampled the posterior of the standard Bayesian lasso with a 201-dimensional Fourier basis and standard improper priors on the training set, and we did the same for the Horseshoe. Next, we sampled the generalized posterior with the learning rate $\hat\eta$ learned by the Safe-Bayesian algorithm, with the same model and priors on the same training set. The grid of $\eta$'s we used was $1, 0.9, 0.8, 0.7, 0.6, 0.5$. We compare the performance of the standard Bayesian lasso and Horseshoe and the Safe-Bayesian versions of the lasso (SB) in terms of mean square error. In all experiments performed with different partitions, priors and number of iterations, SafeBayes never picked $\hat\eta = 1$. We averaged over 10 runs. Moreover, whichever learning rate was chosen by SafeBayes, it always outperformed standard Bayes (with $\eta =1$) in an unchanged set-up. Experiments with different priors for $\lambda$ yielded similar results.


\paragraph{London Air Pollution Data}\label{ap:air}
As training set we use the following data. We start with the first four weeks of the year $2013$, starting at Monday January $7$ at midnight. We have a measurement for (almost) every hour until Sunday February $3^\text{rd}$, $23.00$. We also have data for the first four weeks of $2014$, starting at Monday January $6$ at midnight, until Sunday February $2^\text{nd}$, $23.00$. For each hour in the four weeks we randomly pick a data point from either $2013$ or $2014$. We remove the missing values. We predict for the same time of year in $2015$: starting at Monday January $5$ at midnight, until Sunday February $1^\text{st}$ at $23.00$. We do this with a (Safe-)Bayesian lasso and Horseshoe with a $201$-dimensional Fourier basis and standard improper priors. The grid of $\eta$'s we used for the Safe-Bayesian algorithm was again $1, 0.9, 0.8, 0.7, 0.6, 0.5$. We look at the mean square prediction errors, and average the errors over $20$ runs of the generalized Bayesian lasso with the $\eta$ learned by SafeBayes, and the standard Bayesian lasso and Horseshoe. Again we find that SafeBayes clearly performs better than standard Bayes.

%

\subsection{Details for Figure~\ref{fig:predvar}}
\label{ap:delete}
Here we sampled the posteriors of the standard and generalized Bayesian lasso ($\eta = 0.25$) on $50$ model-wrong data points (approximately half easy points) with $101$ Fourier basis functions, and estimated the predictive variance on a grid of new data points $X_{\text{new}} = \{ -1.00, -0.99, \ldots, 1.00 \}$ with the Monte Carlo estimate:
\begin{equation}
\hat{\textsc{var}}(Y_{\text{new}} \mid X_{\text{new}}, Z_{\text{old}}) = \ex_{\theta \mid Z_{\text{old}}} \left[ \textsc{var}(Y_{\text{new}} \mid \theta) \right] + \hat{\textsc{var}} \left[ \ex(Y_{\text{new}} \mid \theta) \right],
\end{equation}
where
\begin{center}
\begin{gather*}
\ex_{\theta \mid Z_{\text{old}}} \left[ \textsc{var}(Y_{\text{new}} \mid \theta) \right]  = \frac{1}{m} \sum_{k=1}^m \sigma^{2 \left[k \right]}
= \,\overline{\sigma^2},\\ 
\hat{\textsc{var}} \left[ \ex(Y_{\text{new}} \mid \theta) \right] = \,\hat{\textsc{var}} \left[X_{\text{new}}\beta \right]=  \frac{1}{m} \sum_{k=1}^m \left( X_{\text{new}}\beta^{\left[ k \right]} \right)^2 - \left(X_{\text{new}}\overline{\beta}\right)^2.
\end{gather*}
\end{center}
Here $\overline{\beta}$ is the posterior mean of the parameter for the coefficients and $\overline{\sigma^2}$ is the posterior mean of the variance.

\end{document}